\documentclass{article}

\usepackage{authblk}

\title{Goldblatt-Thomason Theorems for Modal Intuitionistic Logics}
\author{Jim de Groot}

\affil{\small
  College of Engineering and Computer Science \\
  The Australian National University \\
  Canberra, ACT, Australia \\
  \texttt{jim.degroot@anu.edu.au}
} 

\date{}

\usepackage{rotate}
\usepackage[dvipsnames]{xcolor}
\usepackage[colorlinks=true,
            linktocpage=true,
            citecolor=RubineRed,
            linkcolor=RoyalBlue]{hyperref}
\usepackage[numbers]{natbib}
\usepackage[a4paper]{geometry}

\usepackage{amsmath}
\usepackage{amssymb}
\usepackage{amsthm}
  \theoremstyle{definition}
  \swapnumbers
    
    \newtheorem{para}{}[section]
    \newtheorem{definition}[para]{Definition}
    \newtheorem{example}[para]{Example}

    \newtheorem{assum}[para]{Assumption}
  \theoremstyle{theorem}
    \newtheorem{lemma}[para]{Lemma}
    
    \newtheorem{theorem}[para]{Theorem}

    \newtheorem{proposition}[para]{Proposition}
    
  \swapnumbers

\usepackage{tikz-cd}
  \usetikzlibrary{arrows,arrows.meta}
  \tikzcdset{arrow style=tikz, diagrams={>=stealth'}}

\usepackage{stmaryrd}
\usepackage{enumerate}
\usepackage{doi}

\DeclareSymbolFont{symbolsC}{U}{txsyc}{m}{n}
\DeclareMathSymbol{\sto}{\mathrel}{symbolsC}{74}
\DeclareMathSymbol{\dright}{\mathrel}{symbolsC}{147}

\usepackage{mathrsfs}
\usepackage[cal=euler,scr=boondoxupr,scrscaled=1.050]{mathalfa}

\newcommand{\mf}[1]{\mathfrak{#1}}
\newcommand{\ms}[1]{\mathscr{#1}}
\newcommand{\fun}[1]{\mathscr{#1}}  % for functors
\newcommand{\cat}[1]{\mathsf{#1}}   % for categories
\newcommand{\mo}[1]{\mathfrak{#1}}   % for models
\newcommand{\amo}[1]{\mathcal{#1}}  % for algebraic models
\newcommand{\lan}[1]{\mathbf{#1}}   % for languages
\renewcommand{\log}[1]{\boldsymbol{\mathcal{#1}}}   % for logics
\newcommand{\fil}[1]{\mathfrak{#1}}
\newcommand{\ff}[1]{\mathfrak{#1}}

\DeclareMathOperator{\id}{id}

\DeclareMathOperator{\Prop}{Prop}
\DeclareMathOperator{\Fr}{Fr}
\DeclareMathOperator{\Alg}{Alg}
\DeclareMathOperator{\Ax}{Ax}
\mathchardef\hyphen="2D
\DeclareMathOperator{\pe}{\mf{pe}_{\tau}}
\newcommand{\power}{\wp}
\renewcommand{\iff}{\quad\text{iff}\quad}
\newcommand{\wiff}{\quad\text{\phantom{iff}}\quad}
\renewcommand{\phi}{\varphi}
\renewcommand{\epsilon}{\varepsilon}

\newcommand{\up}{\fun{u\kern-.1em p}}
\newcommand{\upp}{\fun{u\kern-.1em p'}}

\newcommand{\llb}{\llbracket}
\newcommand{\rrb}{\rrbracket}
\newcommand{\llp}{\llparenthesis\kern.1em}
\newcommand{\rrp}{\kern.1em\rrparenthesis}

\newcommand{\bdlinewidth}{.11ex}

\newcommand{\goodbox}{\hspace{.2ex}\text{%
  \tikz[baseline=-.6ex, rounded corners=.01ex, line width=\bdlinewidth]
    {\draw (-.6ex,-.6ex) rectangle (.6ex,.6ex);}}\kern.2ex}
\newcommand{\gooddiamond}{\hspace{.2ex}\text{%
  \tikz[scale=1.05, baseline=-.6ex, rounded corners=.01ex, rotate=45, line width=\bdlinewidth]
    {\draw (-.5ex,-.5ex) rectangle (.5ex,.5ex);}}\kern.2ex}
\renewcommand{\Box}{\goodbox}
\renewcommand{\Diamond}{\gooddiamond}

\newcommand{\dbox}{\hspace{.2ex}\text{%
  \tikz[baseline=-.6ex, rounded corners=.01ex, line width=\bdlinewidth]
    {\draw (-.6ex,-.6ex) rectangle (.6ex,.6ex);
     \draw[fill=black] (0,0) circle(.12ex);}}\kern.2ex}
\newcommand{\ddiamond}{\hspace{.2ex}\text{%
  \tikz[scale=1.05, baseline=-.6ex, rounded corners=.01ex, rotate=45, line width=\bdlinewidth]
    {\draw (-.5ex,-.5ex) rectangle (.5ex,.5ex);
     \draw[fill=black] (0,0) circle(.11ex);}}\kern.2ex}

\newcommand{\customtri}{\hspace{.2ex}\text{%
  \tikz[baseline=-.6ex, rounded corners=.01ex, line width=.1ex]
    {\draw[-] (0,-.6ex) -- (-.65ex,-.6ex) -- (0,.65ex) -- (.65ex,-.6ex) -- (0,-.6ex);}}\kern.2ex}
\newcommand{\vartri}{{\customtri}}

\newcommand{\dtriangle}{\hspace{.2ex}\text{%
  \tikz[baseline=-.6ex, rounded corners=.01ex, line width=.1ex]
    {\draw[-] (0,-.6ex) -- (-.65ex,-.6ex) -- (0,.65ex) -- (.65ex,-.6ex) -- (0,-.6ex);
     \draw[fill=black,-] (0,-.2ex) circle(.11ex);}}\kern.2ex}

\newcommand{\dheartsuit}{\hspace{-.75ex}\text{%
  \tikz[baseline=-.6ex, rounded corners=.01ex, line width=\bdlinewidth]
    {\draw (0,0) node {$\heartsuit$};
     \draw[fill=black] (0,.15ex) circle(.12ex);}}\kern-.75ex}

\begin{document}

%\begin{frontmatter}
%  \title{Goldblatt-Thomason Theorems for Modal Intuitionistic Logics}
%  \author{Jim de Groot}
%  \address{The Australian National University \\ Canberra \\ Australia}

\maketitle

  \begin{abstract}
    \noindent
    We prove a Goldblatt-Thomason theorem for dialgebraic intuitionistic logics,
    and instantiate it to Goldblatt-Thomason theorems for a wide variety
    of modal intuitionistic logics from the literature.
  \end{abstract}

%  \begin{keyword}
%    Modal logic,
%    intuitionistic logic,
%    Goldblatt-Thomason theorem.
%  \end{keyword}
% \end{frontmatter}

%================================================================================
\section{Introduction}

  A prominent question in the study of (modal) logics and their semantics
  is what classes of frames can be defined as the class of frames satisfying
  some set of formulae.
  Such a class is usually called \emph{axiomatic} or \emph{modally definable}.
  A milestone result partially answering this question in the realm of classical
  normal modal logic is from Goldblatt and Thomason
  and dates back to 1974~\cite{GolTho74}. It states that
  an elementary class of Kripke frames is axiomatic if and only if
  it reflects ultrafilter extensions and is closed under p-morphic images,
  generated subframes and disjoint unions.
  The proof in~\cite{GolTho74} relies on Birkhoff's variety theorem~\cite{Bir35}
  and makes use of the algebraic semantics of the logic.
  A model-theoretic proof was provided almost twenty years later by
  Van Benthem~\cite{Ben93}.

  A similar result for (non-modal) intuitionistic logic was proven by
  Rodenburg \cite{Rod86} (see also \cite{Gol05}), where the interpreting
  structures are \emph{intuitionistic} Kripke frames and models.
  This, of course, requires analogues of the notions of p-morphic images,
  generated subframes, disjoint unions and ultrafilter extensions.
  While the first three carry over straightforwardly from the setting
  of classical normal modal logic,
  ultrafilters need to be replaced by \emph{prime filters}.
  
  In recent years, Goldblatt-Thomason style theorems (which we will simply refer
  to as ``Goldblatt-Thomason theorems'') for many other logics have
  been proven, including for positive normal modal logic \cite{CelJan99},
  graded modal logic \cite{SanMa10},
  modal extensions of {\L}ukasiewicz finitely-valued logics \cite{Teh16},
  LE-logics \cite{ConPalTzi18-arxiv},
  and modal logics with a universal modality \cite{SanVir19}.
  A general Goldblatt-Thomason theorem for coalgebraic logics for
  $\cat{Set}$-coalgebras was given in~\cite{KurRos07}.

  In the present paper we prove Goldblatt-Thomason theorems for
  modal intuitionistic logics.
  These include the extensions of intuitionistic logic with
  a normal modality~\cite{WolZak97,WolZak98,WolZak99},
  a monotone one~\cite[Sec.~6]{Gol93},
  a neighbourhood modality~\cite{DalGreOli20},
  and a strict implication modality~\cite{LitVis18,LitVis19,GroLitPat21}.
  For each we obtain:
  \begin{center}
  {\it
    A class $\ms{K}$ of frames closed under prime filter extensions is axiomatic if \\
    and only if it reflects prime filter extensions and is closed under \\
    disjoint unions, regular subframes and p-morphic images.}
  \end{center}
  Instead of proving each of these results individually, we prove a more
  general Goldblatt-Thomason theorem for \emph{dialgebraic intuitionistic logics},
  merging techniques from~\cite{Gol05} and~\cite{KurRos07}.
  We then apply this to specific instances.
  
  Dialgebraic logic slightly generalises coalgebraic logic and
  was recently introduced in \cite{GroPat20}.
  It provides a framework where modal logics are developed
  parametric in the signature
  of the language and a functor $\fun{T} : \cat{C}' \to \cat{C}$, where
  $\cat{C}'$ is some subcategory of $\cat{C}$.
  While coalgebraic logics are too restrictive to describe modal intuitionistic
  logics (see e.g.~\cite[Rem.~8]{Lit17-arxiv}, \cite[Sec.~2]{GroPat20}),
  the additional flexibility of dialgebraic logic does allow us to
  model a number of them.
  
  The paper is structured as follows.
  In Sec.~\ref{sec:prelim} we recall a semantics for the extension
  of intuitionistic logic with a normal modality $\Box$ from~\cite{WolZak99}.
  Using this as running example, in Sec.~\ref{sec:general} we
  recall the basics of dialgebraic logic and prove the
  Goldblatt-Thomason theorem.
  In particular, this yields a new Goldblatt-Thomason theorem for the
  logic and semantics from Sec.~\ref{sec:prelim}.
  In Sec.~\ref{sec:app} we instantiate the general theorem to several more
  modal intuitionistic
  logics from the literature to obtain new Goldblatt-Thomason theorems.
  
\bigskip\noindent
\textit{Related version.} \,
  This paper is accepted for publications at AiML 2022.

%================================================================================
\section{Normal Modal Intuitionistic Logic}\label{sec:prelim}

  For future reference, we recall
  the extension of
  intuitionistic logic with a unary meet-preserving modality from Wolter and
  Zakharyaschev~\cite{WolZak98,WolZak99}.

\begin{definition}
  Denote the language of intuitionistic logic by $\lan{L}$, with
  proposition letters from some countably infinite set $\Prop$.
  That is, $\lan{L}$ is generated by the grammar
  $$
    \phi ::= \top \mid \bot \mid p \mid \phi \wedge \phi \mid
    \phi \vee \phi \mid \phi \to \phi,
  $$
  where $p \in \Prop$.
  Write $\lan{L}_{\Box}$ for its extension with a unary operator $\Box$.
  Further, let $\log{L}$ be the intuitionistic propositional calculus,
  and let $\log{L}_{\Box}$ be the logic that arises from extending an
  axiomatisation for $\log{L}$ (that we assume includes uniform substitution)
  with the axioms and rule
  \begin{equation}\label{eq:box-ax-rule}
    \Box\top \leftrightarrow \top, \qquad
    \Box p \wedge \Box q \leftrightarrow \Box(p \wedge q), \qquad
    (p \leftrightarrow q)/(\Box p \leftrightarrow \Box q)
  \end{equation}
\end{definition}

  We write $\cat{Pos}$ for the category of posets and order-preserving
  functions. In this paper, we define an \emph{intuitionistic Kripke frame} as a
  poset and we write
  $\cat{Krip}$ for the full subcategory of $\cat{Pos}$ whose morphisms
  are p-morphisms \cite[Sec.~2.1.1]{Bez06}.
  (Sometimes intuitionistic Kripke frames are defined to be preorders.
  For the results presented in this paper there is no discernible difference.)

\begin{definition}
  A \emph{$\Box$-frame} is a triple $(X, \leq, R)$ where $(X, \leq)$ is an
  intuitionistic Kripke frame and $R$ is a relation on $X$ satisfying
  $
    ({\leq} \circ R \circ {\leq}) = R.
  $
  
  Adding a valuation $V : \Prop \to \fun{Up}(X, \leq)$
  ($= \{a \subseteq X \mid x \in a \text{ and } x \leq y \text{ implies } y \in a \}$)
  yields a \emph{$\Box$-model}, in which we can interpret
  $\lan{L}_{\Box}$-formulae.
  Proposition letters are interpreted via the valuation,
  intuitionistic connectives are interpreted as usual in the underlying intuitionistic
  Kripke frame and a state $x$ satisfies $\Box\phi$ if all its
  $R$-successors satisfy $\phi$.
\end{definition}

  While morphisms are not defined in~\cite{WolZak98,WolZak99}, there is
  an obvious choice:

\begin{definition}
  A \emph{$\Box$-morphism} form $(X, \leq, R)$ to $(X', \leq', R')$
  is a function $f : X \to X'$ such that for $E \in \{ \leq, R \}$
  and for all $x, y \in X$ and $z' \in X'$:
  \begin{itemize}
    \item If $xEy$ then $f(x)E'f(y)$;
    \item If $f(x)E'z'$ then $\exists z \in X$ such that $xEz$ and $f(z) = z'$.
  \end{itemize}
  We write $\cat{WZ\Box}$ for the category of $\Box$-frames and -morphisms.
\end{definition}

  The algebraic semantics of $\log{L}_{\Box}$ is given as follows.

\begin{definition}
  A \emph{Heyting algebra with operators} (HAO) is a pair
  $(A, \Box)$ of a Heyting algebra $A$ and a function $\Box : A \to A$
  satisfying $\Box\top = \top$ and $\Box a \wedge \Box b = \Box(a \wedge b)$
  for all $a, b \in A$.
  Together with $\Box$-preserving Heyting homomorphisms,
  these constitute the category $\cat{HAO}$.
\end{definition}
  
  We briefly recall some categories, functors and natural transformations.
  
\begin{definition}\label{def:fun-nat}
  $\cat{DL}$ and $\cat{HA}$ denote the categories of distributive
  lattices and Heyting algebras.
  Let $\up$ be the contravariant functor $\cat{Pos} \to \cat{DL}$
  that sends a poset to the distributive lattice of its upsets and
  an order-preserving function $f$ to $f^{-1}$.
  Write $\fun{pf} : \cat{DL} \to \cat{Pos}$
  for the contravariant functor
  sending $A \in \cat{DL}$ to the set of prime filters of $A$ ordered
  by inclusion, and a homomorphism to its inverse image.
  These restrict to $\upp : \cat{Krip} \to \cat{HA}$
  and $\fun{pf'} : \cat{HA} \to \cat{Krip}$.
  
  Let $\eta : \fun{id}_{\cat{Pos}} \to \fun{pf} \circ \up$
  and $\theta : \fun{id}_{\cat{DL}} \to \up \circ \fun{pf}$
  be the natural trans\-for\-ma\-tions defined by
  $\eta_{(X, \leq)}(x) = \{ a \in \up(X, \leq) \mid x \in a \}$ and
  $\theta_A(a) = \{ \ff{p} \in \fun{pf}A \mid a \in \ff{p} \}$.
  (These are the units of the dual adjunction between $\cat{Pos}$ and
  $\cat{DL}$.)
  Furthermore, $\theta$ restricts to the natural transformation
  $\theta' : \fun{id}_{\cat{HA}} \to \upp \circ \fun{pf'}$.
\end{definition}

  Every $\Box$-frame $(X, \leq, R)$ yields a HAO $(\upp(X, \leq), \Box_R)$
  (called its \emph{complex algebra}),
  with $\Box_R(a) = \{ x \in X \mid xRy \text{ implies } y \in a \}$.
  Conversely, every HAO $(A, \Box)$ gives rise to a $\Box$-frame
  $(\fun{pf'}A, \subseteq, R_{\Box})$, where
  $\ff{p}R_{\Box}\ff{q}$ iff for all $a \in A$, $\Box a \in \ff{p}$ implies $a \in \ff{q}$.
  Concatenating these constructions yields:

\begin{definition}
  The \emph{prime filter extension} of a $\Box$-frame $(X, \leq, R)$
  is the frame $(X^{pe}, \subseteq, R^{pe})$, where $X^{pe}$ is the set of
  prime filters on $(X, \leq)$ and $R^{pe}$ is defined by
  $\ff{p}R^{pe}\ff{q}$ iff for all $a \in \fun{u\kern-.1em p'}(X, \leq)$,
  $\Box_R(a) \in \ff{p}$ implies $a \in \ff{q}$.
\end{definition}

%================================================================================
\section{A General Goldblatt-Thomason Theorem}\label{sec:general}

  We restrict the framework of dialgebraic logic~\cite{GroPat20}
  to an intuitionistic base. Within this, we prove a Goldblatt-Thomason theorem.
  Throughout this section, we show how general constructions specialise to the normal modal
  intuitionistic logic from Sec.~\ref{sec:prelim}.
  Our focus on an intuitionistic propositional base allows us
  to augment the framework of dialgebraic logic from~\cite{GroPat20}
  in the following ways:
  \begin{itemize}
    \item In~\cite{GroPat20} a logic is identified via an initial object
          in some category, which plays the role of the Lindenbaum-Tarski algebra.
          Here we define logics explicitly, by means of an axiomatisation.
    \item Whereas proposition letters in~\cite{GroPat20} are regarded as
          predicate liftings, here we elevate them to a
          special status. This has two reasons: first, it simplifies the
          connection to (frames and models for) modal intuitionistic logics
          from the literature; second, they
          facilitate the use of Birkhoff's variety theorem.
    \item We give dialgebraic definitions of subframes, p-morphic images and
          disjoint unions, and corresponding preservation results.
    \item We give prime filter extensions for models (not just for frames).
  \end{itemize}
  We work towards a Goldblatt-Thomason theorem as follows.
  First we recall the use of dialgebras as frames for modal extensions
  of intuitionistic logic (Sec.~\ref{subsec:lan-frm}),
  and we prove some invariance properties (Sec.~\ref{subsec:invariance}).
  Then we describe algebraic semantics and prime filter extensions
  dialgebraically (Sec.~\ref{subsec:algsem} and~\ref{subsec:pe}).
  This culminates in the Goldblatt-Thomason theorem in Sec.~\ref{subsec:gt}.

%--------------------------------------------------------------------------------
\subsection{Languages and Frames}\label{subsec:lan-frm}

  Dialgebras were introduced by Hagino in~\cite{Hag87} to describe data types.
  Here we use them 
  to describe frames for modal intuitionistic logics.

\begin{definition}
  Let $\fun{F}, \fun{G} : \cat{C} \to \cat{D}$ be functors.
  An \emph{$(\fun{F}, \fun{G})$-dialgebra} is a pair $(X, \gamma)$ where $X \in \cat{C}$
  and $\gamma : \fun{F}X \to \fun{G}X$ is a $\cat{D}$-morphism.
  An \emph{$(\fun{F}, \fun{G})$-dialgebra morphism}
  from $(X, \gamma)$ to $(X', \gamma')$ is a $\cat{C}$-morphism
  $f : X \to X'$ such that $\fun{G}f \circ \gamma = \gamma' \circ \fun{F}f$.
  They constitute the category $\cat{Dialg}(\fun{F}, \fun{G})$.
  In diagrams:
  $$
    \begin{tikzcd}[row sep=-.1em, column sep=1em]
        & \fun{F}X
            \arrow[dd, "\gamma"]
        &
        &
        & \fun{F}X
            \arrow[dd, "\gamma" left]
            \arrow[rr, "\fun{F}f"]
        &
        & \fun{F}X'
            \arrow[dd, "\gamma'"] \\
      \text{objects:}
        &
        &
        & \text{arrows:}
        &
        &
        & \\
        & \fun{G}X
        &
        &
        & \fun{G}X
            \arrow[rr, "\fun{G}f", below]
        &
        & \fun{G}X'
    \end{tikzcd}
  $$
\end{definition}

  We will be concerned with two classes of dialgebras. First,
  $(\fun{i}, \fun{T})$-dialgebras, where $\fun{i} : \cat{Krip} \to \cat{Pos}$
  is the inclusion functor and $\fun{T} : \cat{Krip} \to \cat{Pos}$ is
  any functor, serve as frame semantics for our dialgebraic intuitionistic
  logics.
  Second, dialgebras for functors $\cat{HA} \to \cat{DL}$ 
  will be used as algebraic semantics.
  
\begin{example}
  Let $\fun{P_{up}} : \cat{Krip} \to \cat{Pos}$ be the functor that
  sends an intuitionistic Kripke frame $(X, \leq)$ to its set of
  upsets ordered by reverse inclusion, and a p-morphism
  $f : (X, \leq) \to (X', \leq')$ to
  $\fun{P_{up}}f : \fun{P_{up}}(X, \leq) \to \fun{P_{up}}(X', \leq')
    : a \mapsto f[a]$.
  Then identifying a relation $R$ on $X$ with the map $\gamma_R : (X, \leq) \to \fun{P_{up}}(X, \leq) : x \mapsto \{ y \in X \mid xRy \}$ yields an isomorphism $\cat{WZ\Box} \cong \cat{Dialg}(\fun{i}, \fun{P_{up}})$~\cite[Sec.~2]{GroPat20}.
\end{example}

  Modalities for $\cat{Dialg}(\fun{i}, \fun{T})$ are defined
  via predicate liftings~\cite[Def.~5.7]{GroPat20}.
  
\begin{definition}
  An \emph{$n$-ary predicate lifting} for a functor
  $\fun{T} : \cat{Krip} \to \cat{Pos}$
  is a natural transformation
  $$
    \lambda : (\fun{Up} \circ \fun{i})^n \to \fun{Up} \circ \fun{T}.
  $$
  Here $\fun{Up} : \cat{Pos} \to \cat{Set}$ is the contravariant functor
  that sends a poset to its set of upsets,
  and $(\fun{Up} \circ \fun{i})^n(X, \leq)$ is the $n$-fold product of
  $\fun{Up}(\fun{i}(X, \leq))$ in $\cat{Set}$.
\end{definition}

\begin{definition}
  Let $\Prop$ be a countably infinite set of proposition letters.
  For a set $\Lambda$ of predicate liftings, define the language
  $\lan{L}(\Lambda)$ by the grammar
  $$
    \phi ::= \top
      \mid \bot
      \mid p
      \mid \phi \wedge \phi
      \mid \phi \vee \phi
      \mid \phi \to \phi
      \mid \heartsuit^{\lambda}(\phi_1, \ldots, \phi_n),
  $$
  where $p$ ranges over $\Prop$ and $\lambda \in \Lambda$ is $n$-ary.
\end{definition}
  
\begin{definition}\label{def:interpretation}
  Let $\Lambda$ be a set of predicate liftings for $\fun{T} : \cat{Krip} \to \cat{Pos}$.
  An \emph{$(\fun{i}, \fun{T})$-model} $\mo{M}$ is an
  $(\fun{i}, \fun{T})$-dialgebra $\mo{X} = (X, \leq, \gamma)$
  with a valuation $V : \Prop \to \fun{Up}(X, \leq)$.
  Truth of $\phi \in \lan{L}(\Lambda)$ at $x \in X$ is defined by
  \begin{align*}
    \mo{M}, x \Vdash \top
      &\wiff\text{always} \\
    \mo{M}, x \Vdash \bot
      &\wiff\text{never} \\
    \mo{M}, x \Vdash p
      &\iff x \in V(p) \\
    \mo{M}, x \Vdash \phi \wedge \psi
      &\iff \mo{M}, x \Vdash \phi \text{ and } \mo{M}, x \Vdash \psi \\
    \mo{M}, x \Vdash \phi \vee \psi
      &\iff \mo{M}, x \Vdash \phi \text{ or } \mo{M}, x \Vdash \psi \\
    \mo{M}, x \Vdash \phi \to \psi
      &\iff x \leq y \text{ and } \mo{M}, y \Vdash \phi
                                   \text{ imply } \mo{M}, y \Vdash \psi \\
    \mo{M}, x \Vdash \heartsuit^{\lambda}(\phi_1, \ldots, \phi_n)
      &\iff \gamma(x) \in \lambda_{(X,\leq)}(\llb \phi_1 \rrb^{\mo{M}}, \ldots, \llb \phi_n \rrb^{\mo{M}})
  \end{align*}
  Here $\llb \phi \rrb^{\mo{M}} = \{ x \in X \mid \mo{M}, x \Vdash \phi \}$.
  We write $\mo{M} \Vdash \phi$ if $\mo{M}, x \Vdash \phi$ for all $x \in X$
  and $\mo{X} \Vdash \phi$ if $(\mo{X}, V) \Vdash \phi$ for
  all valuations $V$ for $\mo{X}$.
  If $\Phi \subseteq \lan{L}(\Lambda)$ then we say that
  $\Phi$ is \emph{valid} on $\mo{X}$, and write $\mo{X} \Vdash \Phi$, if
  $\mo{X} \Vdash \phi$ for all $\phi \in \Phi$.
  Also, let
  $$
    \Fr \Phi = \{ \mo{X} \in \cat{Dialg}(\fun{i}, \fun{T}) \mid \mo{X} \Vdash \Phi \}.
  $$
  We call a class $\ms{K} \subseteq \cat{Dialg}(\fun{i}, \fun{T})$
  \emph{axiomatic} if $\ms{K} = \Fr\Phi$ for some $\Phi \subseteq \lan{L}(\Lambda)$.
\end{definition}

\begin{example}\label{exm:pred-lift-box}
  Since $\Box$-frames correspond to $(\fun{i}, \fun{P_{up}})$-dialgebras,
  it is easy to see that $\Box$-models correspond to
  $(\fun{i}, \fun{P_{up}})$-models.
  The modal operator $\Box$ can be induced by the predicate
  lifting $\lambda^{\Box} : \fun{Up} \circ \fun{i} \to \fun{Up} \circ \fun{P_{up}}$
  given by
  $$
    \lambda^{\Box}_{(X, \leq)}
      : \fun{Up}(\fun{i}(X, \leq)) \to \fun{Up}(\fun{P_{up}}(X, \leq))
      : a \mapsto \{ b \in \fun{P_{up}}(X, \leq) \mid b \subseteq a \}.
  $$
  Indeed, if $\mo{M} = (X, \leq, R, V)$ is a $\Box$-model and
  $(X, \leq, \gamma_R, V)$ the corresponding $(\fun{i}, \fun{P_{up}})$-model
  then we have $x \Vdash \Box\phi$ iff every $R$-successor of $x$ satisfies
  $\phi$, i.e.~iff $\gamma_R(x) \subseteq \llb \phi \rrb^{\mo{M}}$.
  By definition the latter is equivalent to
  $\gamma_R(x) \in \lambda^{\Box}_{(X, \leq)}(\llb \phi \rrb^{\mo{M}})$.
\end{example}

  Finally, we define morphisms between $(\fun{i}, \fun{T})$-models.

\begin{definition}
  An \emph{$(\fun{i}, \fun{T})$-model morphism} from $\mo{M} = (\mo{X}, V)$
  to $\mo{M}' = (\mo{X}', V')$ is an $(\fun{i}, \fun{T})$-dialgebra
  morphism $f : \mo{X} \to \mo{X}'$ such that $V = f^{-1} \circ V'$.
\end{definition}

\begin{proposition}\label{prop:mor-pres-truth}
  If $f : \mo{M} \to \mo{M}'$ is an $(\fun{i}, \fun{T})$-model morphism,
  then for all states $x$ of $\mo{M}$ and $\phi \in \lan{L}(\Lambda)$,
  we have $\mo{M}, x \Vdash \phi$ iff $\mo{M}', f(x) \Vdash \phi$.
\end{proposition}
\begin{proof}
  Let $\mo{M} = (X, \leq, \gamma, V)$ and $\mo{M}' = (X', \leq', \gamma', V')$.
  The proof proceeds by induction on the structure of $\phi$.
  If $\phi \in \Prop$ then the claim follows from the
  definition of an $(\fun{i}, \fun{T})$-model morphism.
  The inductive cases for propositional connectives are routine,
  so we focus on the modal case. We restrict our attention to
  unary modalities, higher arities being similar. Compute:
  \begin{align*}
    \mo{M}, x &\Vdash \heartsuit^{\lambda}\phi \\
      &\iff \gamma(x) \in \lambda_{(X, \leq)}(\llb \phi \rrb^{\mo{M}})
            &\text{(Def.~\ref{def:interpretation})} \\
      &\iff \gamma(x) \in \lambda_{(X, \leq)}(f^{-1}(\llb \phi \rrb^{\mo{M}'}))
            &\text{(Induction hypothesis)} \\
      &\iff \gamma(x) \in \lambda_{(X, \leq)}((\fun{i}f)^{-1}(\llb \phi \rrb^{\mo{M}'}))
            &\text{(Because $\fun{i}f = f$)} \\
      &\iff \gamma(x) \in
            (\fun{T}f)^{-1}(\lambda_{(X', \leq')}(\llb \phi \rrb^{\mo{M}'}))
            &\text{(Naturality of $\lambda$)} \\
      &\iff (\fun{T}f)(\gamma(x)) \in \lambda_{(X', \leq')}(\llb \phi \rrb^{\mo{M}'}) \\
      &\iff \gamma'((\fun{i}f)(x)) \in \lambda_{(X', \leq')}(\llb \phi \rrb^{\mo{M}'})
            &\text{($f$ is a dialgebra morphism)} \\
      &\iff \mo{M}', f(x) \Vdash \heartsuit^{\lambda}\phi
            &\text{(Def.~\ref{def:interpretation} and $\fun{i}f = f$)}
  \end{align*}
  This proves the proposition.
\end{proof}

%--------------------------------------------------------------------------------
\subsection{Disjoint Unions, Generated Subframes and p-Morphic Images}
  \label{subsec:invariance}

  The category theoretic analogue of a disjoint union is a coproduct.
  For any $\fun{T} : \cat{Krip} \to \cat{Pos}$
  the category $\cat{Dialg}(\fun{i}, \fun{T})$ has coproducts
  because $\cat{Krip}$ has coproducts and $\fun{i}$ preserves
  them~\cite[Thm.~3.2.1]{Blo12}.
  So we define:

\begin{definition}
  The \emph{disjoint union} of a $K$-indexed family of
  $(\fun{i}, \fun{T})$-dialgebras $\mo{X}_k = (X_k, \leq_k, \gamma_k)$
  is the coproduct $\coprod_{k \in K} \mo{X}_k$ in $\cat{Dialg}(\fun{i}, \fun{T})$.
\end{definition}

\begin{example}
  Let $(X_k, \leq_k, R_k)$ be a $K$-indexed set of $\Box$-frames,
  and $(X_k, \leq_k, \gamma_k)$ the corresponding $(\fun{i}, \fun{P_{up}})$-dialgebras.
  The coproduct $\coprod_{k \in K}(X_k, \leq_k, \gamma_k)$ is given by
  $(X, \leq, \gamma)$, where $(X, \leq)$ is the coproduct of the intuitionistic
  Kripke frames $(X_k, \leq_k)$ (which is computed as in $\cat{Set}$),
  and $\gamma : (X, \leq) \to \fun{P_{up}}(X, \leq)$
  is given by $\gamma(x_k) = \gamma_k(x_k)$ (for $x_k \in X_k$).
  Transforming this back into a $\Box$-frame, we obtain
  $(X, \leq, R)$, with $xRy$ iff there is a $k \in K$ with
  $x, y \in X_k$ and $xR_ky$.
  So this corresponds to the expected notion of disjoint union of
  $\Box$-frames.
\end{example}

\begin{proposition}\label{prop:disj-union}
  Let $\mo{X}_k = (X_k, \leq_k, \gamma_k)$ be a family of
  $(\fun{i}, \fun{T})$-dialgebras indexed by some set $K$.
  Suppose $\mo{X}_k \Vdash \phi$ for all $k \in K$.
  Then $\coprod \mo{X}_k \Vdash \phi$.
\end{proposition}
\begin{proof}
  Let $V$ be a valuation for $\coprod \mo{X}_k$.
  Define the valuation $V_k$ for $\mo{X}_k$ by $V_k(p) = V(p) \cap X_k$.
  Then the coproduct inclusion maps
  $\kappa_k : (\mo{X}_k, V_k) \to (\coprod \mo{X}_k, V)$ are
  $(\fun{i}, \fun{T})$-model morphisms, hence the assumption
  $\mo{X}_k \Vdash \phi$ for all $k \in K$ implies that
  $(\coprod \mo{X}_k, V) \Vdash \phi$.
  Since $V$ was arbitrary, $\coprod \mo{X}_k \Vdash \phi$.
\end{proof}

\begin{definition}\label{def:subf-im}
  Let $\mo{X}' = (X', \leq', \gamma')$ and $\mo{X} = (X, \leq, \gamma)$ be
  $(\fun{i}, \fun{T})$-dialgebras.
  \begin{enumerate}
    \item $\mo{X}'$ is called a \emph{generated subframe} of $\mo{X}$ if there
          exists a p-morphism $f : \mo{X}' \to \mo{X}$ such that
          $f : (X', \leq') \to (X, \leq)$ is an embedding.
    \item $\mo{X}'$ is a \emph{p-morphic image} of $\mo{X}$ if there
          exists a surjective dialgebra morphism $\mo{X} \to \mo{X}'$.
  \end{enumerate}
\end{definition}

\begin{example}
  Guided by \cite[Def.~2.5 and~3.13]{BRV01},
  we could define a generated sub-$\Box$-frame of a $\Box$-frame
  $(X, \leq, R)$ as a $\Box$-frame $(X', \leq', R')$ such that:
  \begin{itemize}
    \item $X' \subseteq X$ and
          ${\leq'} = ({\leq} \cap (X' \times X'))$ and $R' = (R \cap (X' \times X'))$;
    \item if $x \in X'$ and $x \leq y$ or $x R y$, then $y \in X'$.
  \end{itemize}
  With this definition, it can be shown that a $\Box$-frame
  $\mo{X}'$ is isomorphic to a generated sub-$\Box$-frame of a $\Box$-frame $\mo{X}$
  if and only if the dialgebraic rendering of $\mo{X}'$ is a generated
  subframe of the dialgebraic rendering of $\mo{X}$ (as per Def.~\ref{def:subf-im}).
\end{example}

\begin{proposition}\label{prop:subf-im}
  Let $\mo{X}$ be an $(\fun{i}, \fun{T})$-dialgebra such
  that $\mo{X} \Vdash \phi$.
  \begin{enumerate}
    \item If $\mo{X}'$ is a generated subframe of $\mo{X}$
          then $\mo{X}' \Vdash \phi$.
    \item If $\mo{X}'$ is a p-morphic image of $\mo{X}$
          then $\mo{X}' \Vdash \phi$.
  \end{enumerate}
\end{proposition}
\begin{proof}
  We prove the first item, the second item being similar.
  If $\mo{X}' = (X', \leq', \gamma')$ is a generated subframe of
  $\mo{X} = (X, \leq, \gamma)$ then there exists a
  $(\fun{i}, \fun{T})$-dialgebra morphism $f : \mo{X}' \to \mo{X}$
  that is an embedding of the underlying posets.
  Let $V'$ be any valuation for $\mo{X}'$.
  Define a valuation $V^{\uparrow}$ for $\mo{X}$ by
  $V^{\uparrow}(p) = \{ x \in X \mid \exists y \in V'(p) \text{ s.t. } f(y) \leq x \}$.
  Then the fact that $f$ is an embedding implies that $V' = f^{-1} V^{\uparrow}$,
  and therefore $f : (\mo{X}', V') \to (\mo{X}, V^{\uparrow})$ is a
  dialgebra model morphism.
  The assumption that $\mo{X} \Vdash \phi$ together with
  Prop.~\ref{prop:mor-pres-truth} implies that $(\mo{X}', V') \Vdash \phi$.
  Since $V'$ is arbitrary we find $\mo{X}' \Vdash \phi$.
\end{proof}

%--------------------------------------------------------------------------------
\subsection{Axioms and Algebraic Semantics}\label{subsec:algsem}
  
  In order to get intuition for the dialgebraic perspective of algebraic
  semantics, we observe that the category $\cat{HAO}$ is isomorphic to
  a category of dialgebras. In this case, we consider dialgebras for
  functors $\cat{HA} \to \cat{DL}$.
  Again, one of the functors is simply the inclusion functor, which we
  denote by $\fun{j} : \cat{HA} \to \cat{DL}$.

\begin{example}\label{exm:L-box}
  Let $\fun{K} : \cat{HA} \to \cat{DL}$ be the functor that
  sends a Heyting algebra $A$ to the free distributive lattice
  generated by $\{ \dbox a \mid a \in A \}$ modulo $\dbox\top = \top$ and
  $\dbox a \wedge \dbox b = \dbox (a \wedge b)$, where $a$ and $b$ range over $A$.
  The action of $\fun{K}$ on a Heyting homomorphism $h : A \to A'$
  is defined on generators by $\fun{K}h(\dbox a) = \dbox h(a)$.
  Then $\cat{HAO} \cong \cat{Dialg}(\fun{K}, \fun{j})$~\cite[Exm.~3.3]{GroPat20}.
\end{example}
  
  We denote generators by dotted boxes to distinguish
  them from the modality $\Box$.
  Observe that the relations defining $\fun{K}$ correspond
  to the axioms we want a normal box to satisfy.
  We investigate how to generalise this to the setting of some
  arbitrary set $\Lambda$ of predicate liftings for
  a functor $\fun{T} : \cat{Krip} \to \cat{Pos}$.
  
\begin{definition}
  A \emph{rank-1 formula} in $\lan{L}(\Lambda)$ is a formula $\phi$
  such that
  \begin{itemize}
    \item $\phi$ does not contain intuitionistic implication;
    \item each proposition letter appears in the scope of precisely
          one modal operator.
  \end{itemize}
  A \emph{rank-1 axiom} is a formula of the form
  $\phi \leftrightarrow \psi$, where $\phi, \psi$ are rank-1 formulae.
  It is called \emph{sound} if it is valid in all $(\fun{i}, \fun{T})$-dialgebras.
  
  Let $\Ax$ be a collection of sound rank-1 axioms.
  Define the logic $\log{L}(\Lambda, \Ax)$ as the smallest set of
  $\lan{L}(\Lambda)$-formulae containing $\Ax$ and an axiomatisation
  for intuitionistic logic,
  which is closed under modus ponens,
  uniform substitution, and
  $$
    \dfrac{\phi_1 \leftrightarrow \psi_1 \quad \cdots \quad
           \phi_n \leftrightarrow \psi_n}
          {\heartsuit^{\lambda}(\phi_1, \ldots, \phi_n)
           \leftrightarrow \heartsuit^{\lambda}(\psi_1, \ldots, \psi_n)}
           \qquad\text{(congruence rule)}.
  $$
\end{definition}

  Example~\ref{exm:L-box} generalises as follows~\cite[Sec.~5]{GroPat20}.
  
\begin{definition}\label{def:L-lambda-ax}
  Let $\Lambda$ be a set of predicate liftings for $\fun{T}$ and
  $\Ax$ a set of sound rank-1 axioms for $\lan{L}(\Lambda)$.
  For a Heyting algebra $A$, define $\fun{L}^{(\Lambda, \Ax)}A$ to be the
  free distributive lattice generated by
  $\{ \dheartsuit^{\lambda}(a_1, \ldots, a_n) \mid \lambda \in \Lambda, a_i \in A \}$
  modulo the axioms in $\Ax$, where each occurrence of $\heartsuit$ is
  replaced by the formal generator $\dheartsuit$, $\leftrightarrow$ is replaced by $=$,
   and the proposition
  letters range over the elements of $A$.
  (This is well defined since the axioms in $\Ax$ are rank-1 axioms, which
  result in equations constructed from elements of the form
  $\dheartsuit(a_1, \ldots, a_n)$ and distributive lattice connectives.)
  
  If $h : A \to A'$ is a Heyting homomorphism, define
  $\fun{L}^{(\Lambda, \Ax)}h : \fun{L}^{(\Lambda, \Ax)}A \to \fun{L}^{(\Lambda, \Ax)}A'$
  on generators by
  $\fun{L}^{(\Lambda, \Ax)}h(\dheartsuit^{\lambda}(a_1, \ldots, a_n))
    = \dheartsuit^{\lambda}(h(a_1), \ldots, h(a_n))$.
  Then $\fun{L}^{(\Lambda, \Ax)} : \cat{HA} \to \cat{DL}$ defines a
  functor.
\end{definition}

  Again, we use a symbol with a dot in it to denote formal generators,
  and separate them from symbols in the language.

\begin{example}
  Let $\Lambda = \{ \lambda^{\Box} \}$, where $\lambda^{\Box}$ is the
  predicate lifting from Exm.~\ref{exm:pred-lift-box},
  and write $\Box$ instead of $\heartsuit^{\lambda^{\Box}}$.
  Let $\Ax$ consist of the two axioms (not the rule) from \eqref{eq:box-ax-rule},
  and note that these are both rank-1 axioms.
  Then the logic $\log{L}(\Lambda, \Ax)$ coincides with $\log{L}_{\Box}$,
  and the functor obtained from the procedure in Def.~\ref{def:L-lambda-ax}
  is naturally isomorphic to $\fun{K}$ from Exm.~\ref{exm:L-box}.
  (The only difference is the symbol used to represent the formal generators.)
\end{example}

  The following observation allows us to use the Birkhoff variety theorem
  when proving the Goldblatt-Thomason theorem below.

\begin{lemma}\label{lem:dialgLj-variety}
  Let $\fun{L}$ be obtained from predicate liftings and axioms via
  Def.~\ref{def:L-lambda-ax}.
  Then the category $\cat{Dialg}(\fun{L}, \fun{j})$ is a variety of algebras.
\end{lemma}
\begin{proof}
  It is known that the category $\cat{HA}$ of Heyting algebras is a variety
  of algebras. We add to its signature an $n$-ary operation symbol for each
  $n$-ary predicate lifting in $\Lambda$, and to the set of equations defining
  $\cat{HA}$ the equations obtained from $\Ax$ by replacing $\leftrightarrow$ with
  equality and proposition letters with variables.
\end{proof}

  We can evaluate $\lan{L}(\Lambda)$-formulae in a
  $(\fun{L}^{(\Lambda, \Ax)}, \fun{j})$-dialgebra $(A, \alpha)$ with an
  assignment of the proposition letters to elements of $A$.
  Intuitionistic connectives are interpreted as in the Heyting algebra $A$,
  and the interpretation of $\heartsuit^{\lambda}(\phi_1, \ldots, \phi_n)$
  is given by
  $\alpha(\dheartsuit^{\lambda}(\llp \phi_1 \rrp, \ldots, \llp \phi_n \rrp))$,
  where $\llp \phi_i \rrp$ is the interpretation of $\phi_i$.
  We say that $\phi$ is valid in $(A, \alpha)$, and write $(A, \alpha) \models \phi$,
  if $\phi$ evaluates to $\top$ under every assignment of the proposition letters.
  
  This evaluation is closely related to the interpretation of formulae
  in $(\fun{i}, \fun{T})$-dialgebras:
  a formula $\phi$ is valid in some $(\fun{i}, \fun{T})$-dialgebra if and
  only if it is valid in some related algebra, called the complex algebra.

\begin{definition}\label{def:rho}
  Define
  $\rho : \fun{L}^{(\Lambda, \Ax)} \circ \upp
    \to \up \circ \fun{T}$
  on generators by
  $$
    \rho_{(X, \leq)}(\dheartsuit^{\lambda}(a_1, \ldots, a_n))
      = \lambda_{(X, \leq)}(a_1, \ldots, a_n).
  $$
  Then $\rho$ is a well defined transformation because $\Ax$ is assumed to be sound,
  and it is natural because predicate liftings are natural transformations.
  
  It gives rise to a functor
  $(\cdot)^+ : \cat{Dialg}(\fun{i}, \fun{T})
    \to \cat{Dialg}(\fun{L}^{(\Lambda, \Ax)}, \fun{j})$,
  which sends an $(\fun{i}, \fun{T})$-dialgebra $(X, \leq, \gamma)$
  to its \emph{complex algebra} $(\upp(X, \leq), \gamma^+)$, given by
  $$
    \begin{tikzcd}
      \fun{L}^{(\Lambda, \Ax)}(\upp(X, \leq))
            \arrow[r, "\rho_{(X, \leq)}"]
            \arrow[rrr, rounded corners=1.3ex,
                   to path={ -- ([yshift=-2.5ex,xshift=2ex]\tikztostart.south)
                             -- node[above,pos=.56]{\scriptsize$\gamma^+$}
                                ([yshift=-2.5ex,xshift=-1ex]\tikztotarget.south)
                             -- (\tikztotarget.south)}]
        & \up(\fun{T}(X, \leq))
            \arrow[r, "\up\gamma"]
        & [-.5em]
          \up(\fun{i}(X, \leq))
            \arrow[r, equal]
        & [-1.5em]
          \fun{j}(\upp(X, \leq)).
    \end{tikzcd}
  $$
  The action of $(\cdot)^+$ on an $(\fun{i}, \fun{T})$-dialgebra morphism $f$
  is given by $f^+ = \upp f$.
\end{definition}
  
\begin{example}
  Let $(X, \leq, R)$ be a $\Box$-frame and
  $(X, \leq, \gamma)$ the corresponding $(\fun{i}, \fun{P_{up}})$-dialgebra.
  The complex algebra of $(X, \leq, \gamma)$ is the
  $(\fun{K}, \fun{j})$-dialgebra $(\upp(X, \leq), \gamma^+)$,
  where $\gamma^+$ is given by
  $\gamma^+(\dbox a) = \gamma^{-1}(\lambda^{\Box}(a)) = \{ x \in X \mid \gamma(x) \subseteq a \}$.
  Translating this to a HAO, we see that this corresponds precisely
  to the complex algebra of $(X, \leq, R)$ in the sense of
  Sec.~\ref{sec:prelim}.
\end{example}

\begin{proposition}\label{prop:complex-alg}
  Let $\mo{X}$ be an $(\fun{i}, \fun{T})$-dialgebra and
  $\phi \in \lan{L}(\Lambda)$. Then we have
  $$
    \mo{X} \Vdash \phi \iff \mo{X}^+ \models \phi.
  $$
\end{proposition}
\begin{proof}
  This follows from a routine induction on the structure of $\phi$,
  where the base case follows from the fact that valuations for $\mo{X}$
  correspond bijectively to assignments of the proposition letters to
  elements of $\mo{X}^+$.
\end{proof}

%--------------------------------------------------------------------------------
\subsection{Prime Filter Extensions}\label{subsec:pe}

  The proof of the Goldblatt-Thomason theorem relies on Birkhoff's variety
  theorem and the connection between frame semantics and algebraic semantics
  of a logic. As we have seen above, every $\Box$-frame gives rise to a
  complex algebra, or, more generally, every $(\fun{i}, \fun{T})$-dialgebra
  gives rise to a $(\fun{L}, \fun{j})$-dialgebra.
  To transfer the variety theorem from $(\fun{L}, \fun{j})$-dialgebras back to
  $(\fun{i}, \fun{T})$-dialgebras, we need a functor
  $(\cdot)_+ : \cat{Dialg}(\fun{L}, \fun{j}) \to \cat{Dialg}(\fun{i}, \fun{T})$
  such that
  for each $(\fun{i}, \fun{T})$-dialgebra $\mo{X}$,
  \begin{equation}\label{eq:goal2}\tag{$\star$}
    (\mo{X}^+)_+ \Vdash \phi
    \quad\text{implies}\quad\mo{X} \Vdash \phi.
  \end{equation}

\begin{assum}\label{ass:3.4}
  Throughout this subsection, let
  $\fun{T} : \cat{Krip} \to \cat{Pos}$ be a functor,
  $\Lambda$ a set of predicate liftings for $\fun{T}$,
  and a set $\Ax$ of sound rank-1 axioms from $\lan{L}(\Lambda)$.
  Abbreviate $\fun{L} := \fun{L}^{(\Lambda, \Ax)}$ and
  $\rho := \rho^{(\Lambda, \Ax)}$.
\end{assum}

  A functor $(\cdot)_+ : \cat{Dialg}(\fun{L}, \fun{j}) \to \cat{Dialg}(\fun{i}, \fun{T})$
  arises from a natural transformation
  $\tau$ in the same way as $\rho$ induced a functor from frames to
  complex algebras. To stress its dependence on the choice of $\tau$,
  we denote it by $(\cdot)_{\tau}$ instead of $(\cdot)_+$.
  
\begin{definition}\label{def:tau}
  Let $\tau : \fun{pf} \circ \fun{L} \to \fun{T} \circ \fun{pf'}$
  be a natural transformation.
  Then we define the contravariant functor
  $(\cdot)_{\tau} : \cat{Dialg}(\fun{L}, \fun{j}) \to \cat{Dialg}(\fun{i}, \fun{T})$
  on objects by sending a $(\fun{L}, \fun{j})$-dialgebra $\amo{H} = (H, \alpha)$ to
  the $(\fun{i}, \fun{T})$-dialgebra $\amo{H}_{\tau}$ given by
  $$
    \begin{tikzcd}
      \fun{i}(\fun{pf'}H)
            \arrow[r, equal]
        & [-1.5em]
          \fun{pf}(\fun{j}H)
            \arrow[r, "\fun{pf}\alpha"]
        & \fun{pf}(\fun{L}H)
            \arrow[r, "\tau_H"]
        & \fun{T}(\fun{pf'}H).
    \end{tikzcd}
  $$
  For a $(\fun{L}, \fun{j})$-dialgebra morphism $h : \amo{H} \to \amo{H}'$
  we define $h_{\tau} = \fun{pf'}h : \amo{H}'_{\tau} \to \amo{H}_{\tau}$.
  Naturality of $\tau$ ensures that this is well defined.
\end{definition}

  We call $(\mo{X}^+)_{\tau}$ the $\tau$-prime filter extension of
  an $(\fun{i}, \fun{T})$-dialgebra $\mo{X}$ if $\tau$ satisfies a
  sufficient condition that ensures that~\eqref{eq:goal2} holds
  (by Prop.~\ref{prop:pf-key}).
  This condition relies on the following variation of the adjoint mate of $\rho$.

\begin{definition}
  Let $\rho : \fun{L} \circ \upp \to \up \circ \fun{T}$.
  Then we write $\rho^{\flat}$ for the natural transformation defined as
  the composition
  $$
    \begin{tikzcd}
      \fun{T} \circ \fun{pf'}
            \arrow[r, "\eta_{\fun{T} \circ \fun{pf'}}"]
        & \fun{pf} \circ \up \circ \fun{T} \circ \fun{pf'}
            \arrow[r, "\fun{pf}\rho_{\fun{pf'}}"]
        & \fun{pf} \circ \fun{L} \circ \upp \circ \fun{pf'}
            \arrow[r, "\fun{pf}(\fun{L}\theta')"]
        & \fun{pf} \circ \fun{L},
    \end{tikzcd}
  $$
  where $\eta$ and $\theta$ are defined as in Def.~\ref{def:fun-nat}.
\end{definition}

\begin{definition}\label{def:pfe5}
  Let $\tau$ be a natural transformation such that
  $\rho^{\flat} \circ \tau = \fun{id}_{\fun{pf} \circ \fun{L}}$.
  \begin{enumerate}
    \item Define
          $\pe := (\cdot)_{\tau} \circ (\cdot)^+
            : \cat{Dialg}(\fun{i}, \fun{T})\to \cat{Dialg}(\fun{i}, \fun{T})$.
          We call $\pe\mo{X}$ the \emph{$\tau$-prime filter extension}
          of $\mo{X} \in \cat{Dialg}(\fun{i}, \fun{T})$.
    \item The \emph{$\tau$-prime filter extension} of a model
          $\mo{M} = (\mo{X}, V)$ is $\pe\mo{M} := (\pe\mo{X}, V^{pe})$,
          where $V^{pe}(p) = \{ \ff{q} \in \pe\mo{X} \mid V(p) \in \ff{q} \}$
          for all $p \in \Prop$.
  \end{enumerate}
\end{definition}

  Observe that the prime filter extension of an $(\fun{i}, \fun{T})$-dialgebra
  $\mo{X} = (X, \leq, \gamma)$ is of the form
  $\pe\mo{X} = (X^{pe}, \subseteq, \gamma^{pe})$,
  where $X^{pe}$ denotes the set of prime filters of upsets of
  $(X, \leq)$ and $\gamma^{pe}$ is
  computed using both $\rho$ and $\tau$.
  
  We now show that $\tau$-prime filter extensions satisfy~\eqref{eq:goal2}.

\begin{proposition}\label{prop:pf-key}
  Let $\tau$ be a natural transformation such that
  $\rho^{\flat} \circ \tau = \fun{id}_{\fun{pf} \circ \fun{L}}$,
  $\mo{X} = (X, \leq, \gamma)$ an $(\fun{i}, \fun{T})$-dialgebra, $\mo{M} = (\mo{X}, V)$
  a model based on $\mo{X}$, $\phi \in \lan{L}(\Lambda)$.
  %Then:
  \begin{enumerate}
    \item For all prime filters $\ff{q} \in X^{pe}$ we have
          $\pe\mo{M}, \ff{q} \Vdash \phi$ iff $\llb \phi \rrb^{\mo{M}} \in \ff{q}$.
    \item For all states $x \in X$ we have
          $\mo{M}, x \Vdash \phi$ iff $\pe\mo{M}, \eta_{(X, \leq)}(x) \Vdash \phi$.
    \item If $\pe\mo{X} \Vdash \phi$ then $\mo{X} \Vdash \phi$.
  \end{enumerate}
\end{proposition}
\begin{proof}
  The proof of the proposition is given in the appendix.
\end{proof}

\begin{example}\label{exm:box-frame-tau}
  Returning to our example of $\Box$-frames, we wish to find a natural
  transformation $\tau^{\Box}$ such that
  $(\rho^{\Box})^{\flat} \circ \tau^{\Box} = \fun{id}_{\fun{pf} \circ \fun{L}^{\Box}}$.

  Before defining $\tau^{\Box}$, let us get an idea of what
  $(\rho^{\Box})^{\flat}$ looks like.
  Let $A$ be a Heyting algebra and $Q \in \fun{pf}(\fun{L}^{\Box}A)$.
  Since $Q$ is determined by elements of the form $\dbox a$ it contains,
  where $a \in A$,
  we pay special attention to these elements.
  For $D \in \fun{P_{up}} \circ \fun{pf'}A$ and $a \in A$ we have
  \begin{align*}
    \dbox a \in (\rho^{\Box}_A)^{\flat}(D)
      &\iff \rho_{\fun{pf'}A}((\fun{pf}(\fun{L}^{\Box}\theta'_A))(\dbox a)) \in \eta_{\fun{T}(\fun{pf'}A)}(D) \\
      &\iff \rho_{\fun{pf'}A}(\dbox\theta'_A(a)) \in \eta_{\fun{T}(\fun{pf'}A)}(D) \\
      &\iff D \in \rho_{\fun{pf'}A}(\dbox \theta'_A(a)) \\
      &\iff D \subseteq \theta'_A(a)
  \end{align*}
  Guided by this we define
  $\tau : \fun{pf} \circ \fun{L}^{\Box} \to \fun{P_{up}} \circ \fun{pf'}$
  on components by
  $$
    \tau_A^{\Box}
      : \fun{pf}(\fun{L}^{\Box}A) \to \fun{P_{up}}(\fun{pf'}A)
      : Q \mapsto \{ \ff{p} \in \fun{pf'}A \mid \forall a \in A, \dbox a \in Q \text{ implies } a \in \ff{p} \}
  $$
  With this definition we can prove the following lemma,
  the proof of which can be found in the appendix.
  
\begin{lemma}\label{lem:Box-canonical}
  $\tau^{\Box}$ is a natural transformation such that
  $(\rho^{\Box})^{\flat} \circ \tau^{\Box} = \fun{id}_{\fun{pf} \circ \fun{L}^{\Box}}$.
\end{lemma}

  Now suppose $(A, \Box)$ is a HAO, and $\amo{A} = (A, \alpha)$ its corresponding
  $(\fun{L}^{\Box}, \fun{j})$-dialgebra
  (with $\alpha$ given by $\alpha(\dbox a) = \Box a$).
  We have $\amo{A}_{\tau} = (\fun{pf'}A, \subseteq, \gamma)$, where
  $$
    \gamma(\ff{q})
      = \{ \ff{p} \in \fun{pf'}A \mid \forall a \in A,
           \dbox a \in \alpha^{-1}(\ff{q}) \text{ implies } a \in \ff{p} \}.
  $$
  Note that $\dbox a \in \alpha^{-1}(\ff{q})$ iff $\Box a = \alpha(\dbox a) \in \ff{q}$.
  Therefore, translating $\gamma$ to a relation $R_{\gamma}$, we obtain:
  $\ff{q}R_{\gamma}\ff{p}$ iff $\Box a \in \ff{q}$ implies $a \in \ff{p}$ for
  all $a \in A$.
  
  It follows that the $(\fun{i}, \fun{T})$-dialgebra
  corresponding to the prime filter extension of a $\Box$-frame $(X, \leq, R)$
  (as in Sec.~\ref{sec:prelim})
  coincides with the $\tau^{\Box}$-prime filter extension of the dialgebraic
  rendering of $\mo{X}$.
  So, modulo dialgebraic translation, prime filter extensions and 
  $\tau^{\Box}$-prime filter extensions of $\Box$-frames coincide.
\end{example}

%--------------------------------------------------------------------------------
\subsection{The Goldblatt-Thomason Theorem}\label{subsec:gt}

  Finally, we put our theory to work and prove a Goldblatt-Thomason theorem
  for dialgebraic intuitionistic logics.
  We work with the same assumptions as in Assum.~\ref{ass:3.4}.
  Additionally, we assume that we have a natural transformation
  $\tau : \fun{pf} \circ \fun{L} \to \fun{T} \circ \fun{pf'}$
  such that $\rho^{\flat} \circ \tau = \fun{id}_{\fun{pf} \circ \fun{L}}$.
  This allows us to use Def.~\ref{def:pfe5}.

\begin{definition}
  If $\Phi \subseteq \lan{L}(\Lambda)$ and $\amo{A} \in \cat{Dialg}(\fun{L}, \fun{j})$
  then we write $\amo{A} \models \Phi$ if $\amo{A} \models \phi$ for all
  $\phi \in \Phi$. Besides, we let
  $
    \Alg \Phi = \{ \amo{A} \in \cat{Dialg}(\fun{L}, \fun{j}) \mid \amo{A} \models \Phi \}
  $
  be the collection of $(\fun{L}, \fun{j})$-dialgebras satisfying $\Phi$.
  We say that a class $\ms{C} \subseteq \cat{Dialg}(\fun{L}, \fun{j})$
  is \emph{axiomatic} if $\ms{C} = \Alg \Phi$
  for some collection $\Phi$ of $\lan{L}(\Lambda)$-formulae.
\end{definition}

\begin{lemma}\label{lem:alg-axiomatic-variety}
  $\ms{C} \subseteq \cat{Dialg}(\fun{L}, \fun{j})$ is axiomatic
  iff it is a variety of algebras.
\end{lemma}
\begin{proof}
  If
  $\cat{A} = \{ \amo{A} \in \cat{Dialg}(\fun{L}, \fun{j}) \mid \amo{A} \models \Phi \}$, 
  then it is precisely the variety of algebras satisfying
  $\phi^x \leftrightarrow \top$, where $\phi \in \Phi$ and $\phi^x$ is
  the formula we get from $\phi$ by replacing the proposition letters with
  variables from some set $S$ of variables.
  Conversely, suppose $\cat{A}$ is a variety of algebras given by a set
  $E$ of equations using variables in $S$. For each equation $\phi = \psi$
  in $E$, let $(\phi \leftrightarrow \psi)^p$ be the formula we get from
  replacing the variables in $\phi \leftrightarrow \psi$ with proposition
  letters. Then we have
  $\cat{A} = \Alg \{ (\phi \leftrightarrow \psi)^p \mid \phi = \psi \in E \}$.
\end{proof}

  For a class $\ms{K}$ of $(\fun{i}, \fun{T})$-dialgebras, write
  $\ms{K}^+ = \{ \mo{X}^+ \mid \mo{X} \in \ms{K} \}$ for the collection
  of corresponding complex algebras.
  Also, if $\ms{C}$ is a class of algebras, then we write
  $H\ms{C}$, $S\ms{C}$ and $P\ms{C}$ for its closure under
  \textbf{h}omomorphic images, \textbf{s}ubalgebras and \textbf{p}roducts,
  respectively.

\begin{lemma}\label{lem:mdv}
  A class $\ms{K} \subseteq \cat{Dialg}(\fun{i}, \fun{T})$ is axiomatic
  if and only if
  \begin{equation}\label{eq:mdv}
    \ms{K} = \{ \mo{X} \in \cat{Dialg}(\fun{i}, \fun{T})
                \mid \mo{X}^+ \in HSP(\ms{K}^+) \}.
  \end{equation}
\end{lemma}
\begin{proof}
  Suppose $\ms{K}$ is axiomatic, i.e.~$\ms{K} = \Fr\Phi$.
  Then it follows from Prop.~\ref{prop:complex-alg} and the fact that $H$,
  $S$ and $P$ preserve validity of formulae that
  \eqref{eq:mdv} holds.
  Conversely, suppose \eqref{eq:mdv} holds.
  Since $HSP(\ms{K}^+)$ is a variety, Birkhoff's variety theorem
  states that it is of the from $\Alg\Phi$.
  It follows that $\ms{K} = \Fr\Phi$.
\end{proof}

  We now have all the ingredients to prove the Goldblatt-Thomason theorem.

\begin{theorem}\label{thm:gt}
  Let $\ms{K} \subseteq \cat{Dialg}(\fun{i}, \fun{T})$ be closed under 
  $\tau$-prime filter extensions.
  Then $\ms{K}$ is axiomatic if and only if $\ms{K}$ reflects $\tau$-prime
  filter extensions and is closed under disjoint
  unions, generated subframes and p-morphic images.
\end{theorem}
\begin{proof}
  The implication from left to right follows from
  Sec.~\ref{subsec:invariance} and Prop.~\ref{prop:pf-key}.
  For the converse, by Lem.~\ref{lem:mdv} it suffices to prove that
  $\ms{K} = \{ \mo{X} \in \cat{Dialg}(\fun{i}, \fun{T}) \mid \mo{X}^+ \in HSP(\ms{K}^+) \}$.
  So let $\mo{X} = (X, \gamma) \in \cat{Dialg}(\fun{i}, \fun{T})$ and
  suppose $\mo{X}^+ \in HSP(\ms{K}^+)$.
  Then there are $\mo{Z}_i \in \ms{K}$ such that $\mo{X}^+$ is the
  homomorphic image of a sub-dialgebra $\amo{A}$ of the
  product of the $\mo{Z}_i^+$.
  In a diagram:
  $$
    \begin{tikzcd}
      \mo{X}^+
        & [2.5em]
          \amo{A}
            \arrow[l, ->>, "\text{surjective}" {above,pos=.46}]
            \arrow[r, >->, "\text{injective}"]
        & [2.5em]
          \prod \mo{Z}_i^+
    \end{tikzcd}
  $$
  Since $\prod\mo{Z}_i^+ = (\coprod \mo{Z}_i)^+$,
  dually this yields
  $$
    \begin{tikzcd}
      (\mo{X}^+)_{\tau}
            \arrow[r, >->, "\text{gen.~subframe}"]
        & [5em]
          \amo{A}_{\tau}
        & [5em]
          \big(\big(\coprod \mo{Z}_i\big)^+\big)_{\tau}
            \arrow[l, ->>, "\text{p-morphic image}" {above,pos=.46}]
    \end{tikzcd}
  $$
  We have $\coprod \mo{Z}_i \in \ms{K}$ because $\ms{K}$ is closed under coproducts,
  and $\big((\coprod \mo{Z}_i)^+\big)_{\tau} \in \ms{K}$ because $\ms{K}$ is
  closed under prime filter extensions.
  Then $\amo{A}_{\tau} \in \ms{K}$ and $(\mo{X}^+)_{\tau} \in \ms{K}$ because 
  $\ms{K}$ is closed under p-morphic images and generated subframes.
  Finally, since $\ms{K}$ reflects prime filter extensions we find
  $\mo{X} \in \ms{K}$.
\end{proof}

  Circling back to $\Box$-frames, it follows from Lem.~\ref{lem:Box-canonical}
  and Thm.~\ref{thm:gt} that:

\begin{theorem}\label{thm:gt-frm}
  Suppose $\ms{K} \subseteq \cat{WZ\Box}$ is closed under prime filter extensions.
  Then $\ms{K}$ is axiomatic if and only if it reflects prime filter extensions
  and is closed under disjoint unions, generated subframes and p-morphic
  images.
\end{theorem}

%================================================================================
\section{Applications}\label{sec:app}

  In each of the following subsection we recall a modal intuitionistic logic and
  model it dialgebraically.
  We use this to derive a notion of prime filter extension and
  we apply Thm.~\ref{thm:gt} to obtain a Goldblatt-Thomason theorem.

%--------------------------------------------------------------------------------
\subsection{Goldblatt's Geometric Modality I}\label{subsec:mon-I}

  The extension of intuitionistic logic with a monotone modality,
  here denoted by $\vartri$, was first studied by Goldblatt
  in~\cite[Sec.~6]{Gol93}. It is closely related to its classical
  counterpart~\cite{Che80,Han03,HanKup04}, except that the underlying
  propositional logic is intuitionistic.
  A dialgebraic perspective was given in~\cite[Sec.~8]{GroPat20}.
  
  Let $\lan{L}_{\vartri}$ denote the language of intuitionistic logic
  extended with a unary operator $\vartri$, and write
  $\log{L}_{\vartri}$ for the logic obtained from extending
  intuitionistic logic $\log{L}$ with the axiom
  $
    \vartri(p \wedge q) \to \vartri p
  $
  and the congruence rule for $\vartri$.

\begin{definition}
  An \emph{intuitionistic monotone frame} (or \emph{IM-frame})
  is a triple $(X, \leq, N)$ where
  $(X, \leq)$ is an intuitionistic Kripke frame and $N$ is a function that
  assigns to each $x \in X$ a collection of upsets of $(X, \leq)$ such that:
  \begin{itemize}
    \item if $a \in N(x)$ and $a \subseteq b \in \fun{Up}(X, \leq)$,
          then $b \in N(x)$;
    \item if $x \leq y$ then $N(x) \subseteq N(y)$.
  \end{itemize}
  An \emph{intuitionistic monotone frame morphism} (IMF-morphism) from
  $(X_1, \leq_1, N_1)$ to $(X_2, \leq_2, N_2)$ is a p-morphism
  $f : (X_1, \leq_1) \to (X_2, \leq_2)$ such that
  $f^{-1}(a_2) \in N_1(x_1)$ iff $a_2 \in N_2(f(x_1))$
  for all $x_1 \in X_1$ and $a_2 \in \fun{Up}(X_2, \leq_2)$.
  We write $\cat{Mon}$ for the category of intuitionistic monotone
  frames and morphisms.
\end{definition}

  An \emph{intuitionistic monotone model} is a tuple $\mo{M} = (X, \leq, N, V)$
  such that $(X, \leq, N)$ is an intuitionistic monotone frame and
  $V : \Prop \to \fun{Up}(X, \leq)$ is a valuation.
  The interpretation of $\lan{L}_{\vartri}$-formulae at a state $x$
  in $\mo{M}$ is defined recursively, where the propositional cases
  are as usual and
  $\mo{M}, x \Vdash \vartri \phi$ iff $\llb \phi \rrb^{\mo{M}} \in N(x)$.
  We now take a dialgebraic perspective.

\begin{definition}
  For an intuitionistic Kripke frame $(X, \leq)$, define
  $$
    \fun{M}(X, \leq) = 
    \{ W \subseteq \fun{Up}(X, \leq) \mid \text{ if } a \in W \text{ and }
       a \subseteq b \in \fun{Up}(X, \leq) \text{ then } b \in W \}
  $$
  ordered by inclusion.
  For a p-morphism $f : (X_1, \leq_1) \to (X_2, \leq_2)$, let
  $$
    \fun{M}f
      : \fun{M}(X_1, \leq_1) \to \fun{M}(X_2, \leq_2)
      : W \mapsto \{ a_2 \in \fun{Up}(X_2, \leq_2) \mid f^{-1}(a_2) \in W \}.
  $$
  Then $\fun{M} : \cat{Krip} \to \cat{Pos}$ defines a functor.
\end{definition}

\begin{theorem}[\cite{GroPat20}, Thm.~8.3]\label{thm:Mon-dialg}
  We have $\cat{Mon} \cong \cat{Dialg}(\fun{i}, \fun{M})$.
\end{theorem}

  Translating the dialgebraic notion of disjoint union to IM-frames gives:

\begin{definition}\label{def:mif-coprod}
  Let $\{ (X_k, \leq_k, N_k) \mid k \in K \}$ be a $K$-indexed set
  of IM-frames.
  The disjoint union $\coprod_{k \in K}(X_k, \leq_k, N_k)$
  is the frame $(X, \leq, N)$ where $(X, \leq)$ is the disjoint union of the
  intuitionistic Kripke frames $(X_k, \leq_k)$, and
  $N$ is given by $a \in N(x_k)$ iff $a \cap X_k \in N_k(x_k)$ for all
  $a \in \fun{Up}(X, \leq)$ and $x_k \in X_k$.  
\end{definition}

\begin{definition}
  An IM-frame $\mo{X}'$ is a \emph{generated subframe} of an IM-frame
  $\mo{X}$ if there exists
  an IMF-morphism $\mo{X}' \to \mo{X}$ that is an embedding of posets,
  and $\mo{X}'$ is a \emph{p-morphic image} of $\mo{X}$ if there
  is a surjective IMF-morphism $\mo{X} \to \mo{X}'$.
\end{definition}

  The modal operator
  $\vartri$ can be introduced by the predicate lifting
  $\lambda^{\vartri} : \fun{Up} \circ \fun{i} \to \fun{Up} \circ \fun{M}$ given by
  $$
    \lambda^{\vartri}_{(X, \leq)}(a) = \{ W \in \fun{M}(X, \leq) \mid a \in W \}.
  $$
  With $\Ax = \{ \vartri(a \wedge b) \wedge \vartri a \leftrightarrow \vartri(a \wedge b) \}$ we have
  $\log{L}_{\vartri} = \log{L}(\{ \lambda^{\vartri} \}, \Ax)$.
  Its algebraic semantics is given by $(\fun{L}^{\vartri}, \fun{j})$-dialgebras,
  where $\fun{L}^{\vartri} : \cat{HA} \to \cat{DL}$ is the functor sending
  $A$ to the free distributive lattice generated by
  $\{ \dtriangle a \mid a \in H \}$ modulo $\dtriangle(a \wedge b) \leq \dtriangle b$.
  The corresponding natural transformation
  $\rho^{\vartri} : \fun{L}^{\vartri} \circ \upp \to \up \circ \fun{M}$
  is defined on generators
  by $\rho^{\vartri}_{(X, \leq)}(a) = \{ W \in \fun{M}(X, \leq) \mid a \in W \}$.
  
  Towards prime filter extensions and a Goldblatt-Thomason
  theorem we need to define a right inverse $\tau$ of $(\rho^{\vartri})^{\flat}$.
  To garner inspiration we investigate what
  $(\rho^{\vartri})^{\flat}_A : \fun{M}(\fun{pf'}A) \to \fun{pf}(\fun{L}^{\vartri}A)$
  looks like for $A \in \cat{HA}$.
  We have
  \begin{align*}
    \dtriangle a \in (\rho^{\vartri})^{\flat}_A(W)
       \iff \rho_{\fun{pf'}A}(\dtriangle \theta_A'(a))
           \in \eta_{\fun{M} \circ \fun{pf'}A}(W)
       \iff \theta_A'(a) \in W
  \end{align*}
  for all $W \in \fun{M}(\fun{pf'}A)$ and $a \in A$.
  (Recall that $\theta'_A(a) = \{ \ff{q} \in \fun{pf'}A \mid a \in \ff{q} \}$.)

\begin{definition}\label{def:closed-open}
  Let $A \in \cat{HA}$. We call $D \in \upp(\fun{pf'}A)$
  \emph{closed} if $D = \bigcap \{ \theta_A'(a) \mid a \in A \text{ and } D \subseteq \theta_A'(a) \}$, and
  \emph{open} if
  $D = \bigcup \{ \theta_A'(a) \mid a \in A \text{ and } \theta_A'(a) \subseteq D \}$.
\end{definition}

  (Indeed, this coincides with closed and open upsets of $\fun{pf'}A$, conceived
  of as an Esakia space~\cite[Sec.~2.3.3]{Bez06}.)
  Upsets of the form $\theta'_A(a)$ are closed \emph{and} open.
  
\begin{definition}\label{def:mon-tau}
  For a Heyting algebra $A$, define
  $\tau_A : \fun{pf} \circ \fun{L}^{\vartri}A \to \fun{M} \circ \fun{pf'}A$
  as follows. Let $Q \in \fun{pf} ( \fun{L}^{\vartri}A)$ and
  $D \in \fun{Up}(\fun{pf'}A)$, and define:
  \begin{itemize}
    \item If $D = \theta'_A(a)$ for some
          $a \in A$, then $\theta'_A(a) \in \tau_A(Q)$ if $\dtriangle  a \in Q$;
    \item If $D$ is closed
          then $D \in \tau_A(Q)$ if for all $a \in A$, $D \subseteq \theta'_A(a)$ implies
          $\dtriangle a \in Q$.
    \item For other $D$, $D \in \tau_A(Q)$ if there is a closed upset
          $C \subseteq D$ such that $C \in \tau_A(Q)$.
  \end{itemize}
\end{definition}

  It is easy to see that $\tau^{\vartri}_A$ is an order-preserving function,
  i.e.~a morphism in $\cat{Pos}$.
  The next lemma states that $\tau^{\vartri}$ is a natural transformation.
  We postpone the unexciting proof to the appendix.

\begin{lemma}\label{lem:mon-tau-natural}
  The transformation $\tau$ from Def.~\ref{def:mon-tau}
  is natural.
  Moreover, $(\rho^{\vartri})^{\flat}_A \circ \tau^{\vartri}_A = \fun{id}(\fun{pf}(\fun{L}^{\vartri}A))$ for every Heyting algebra $A$.
\end{lemma}

  Translating the dialgebraic definition of a prime filter extension
  to IM-frames gives a definition of prime filter extension for
  IM-frames. We emphasise that this definition relies on $\tau^{\vartri}$.
  In the next section we derive a different notion of prime filter extension
  for IM-frames, with its own Goldblatt-Thomason theorem.

\begin{definition}
  The \emph{$\tau^{\vartri}$-prime filter extension} of an IM-frame
  $(X, \leq, N)$ is the IM-frame $(X^{pe}, \subseteq, N^{pe})$, where
  $N^{pe}$ is given as follows.
  Let $\vartri_N(a) = \{ x \in X \mid a \in N(x) \}$,
  and for $\ff{q} \in X^{pe}$ and $D \in \fun{Up}(X^{pe}, \subseteq)$
  define:
  \begin{itemize}
    \item If $a \in \upp(X, \leq)$,
          then $\theta_A'(a) \in N^{pe}(\ff{q})$ if $\vartri_Na \in \ff{q}$;
    \item If $D$ is closed then 
          $D \in N^{pe}(\ff{q})$ if $\theta_A'(a) \in N^{pe}(\mf{q})$
          for all $\theta_A'(a)$ containing $D$;
    \item For any $D$, $D \in N^{pe}(\ff{q})$ if there is a closed
          $C \subseteq D$ such that $C \in N^{pe}(\ff{q})$.
  \end{itemize}
\end{definition}

  Now Thm.~\ref{thm:gt} instantiates to:

\begin{theorem}
  Suppose $\ms{K}$ is a class of IM-frames closed under $\tau^{\vartri}$-prime filter extensions.
  Then $\ms{K}$ is axiomatic iff it reflects $\tau^{\vartri}$-prime filter extensions
  and is closed under disjoint unions, generated subframes and p-morphic
  images.
\end{theorem}

%--------------------------------------------------------------------------------
\subsection{Goldblatt's Geometric Modality II}

  We substantiate the claim that a logic may have several notions of
  prime filter extension by giving a different
  right-inverse of $(\rho^{\vartri})^{\flat}$
  from Sec.~\ref{subsec:mon-I}. The setup
  is the same as in Sec.~\ref{subsec:mon-I}, so we
  proceed by defining a right-inverse of $(\rho^{\vartri})^{\flat}$.

\begin{definition}
  For a Heyting algebra $A$, define
  $\sigma_A : \fun{pf} \circ \fun{L}^{\vartri}A \to \fun{M} \circ \fun{pf'}A$
  by sending $Q \in \fun{pf}(\fun{L}^{\vartri}A)$ to $\sigma_A(Q)$, where:
  \begin{itemize}
    \item For open upsets $D$,
          let $D \in \sigma_A(Q)$ if $\exists a \in A$ s.t.~$\dtriangle a \in Q$
          and $\theta_A'(a) \subseteq D$;
    \item For any other upset $D$, let $D \in \sigma_A(Q)$ if all open supersets
          of $D$ are in $\sigma_A(Q)$.
  \end{itemize}
\end{definition}

  Similar to Lem.~\ref{lem:mon-tau-natural} we can prove the following.

\begin{lemma}\label{lem:sigma-natural}
  $\sigma = (\sigma_A)_{A \in \cat{HA}} : \fun{pf} \circ \fun{L}^{\vartri} \to \fun{M} \circ \fun{pf'}$
  is a natural transformation,
  and for every Heyting algebra $A$, we have
  $\rho^{\flat}_A \circ \sigma_A = \fun{id}_{\fun{pf}(\fun{L}^{\vartri}A)}$.
\end{lemma}

  Now $\sigma$ yields a different notion of
  prime filter extension, the precise definition of which we leave to the reader.
  Thm.~\ref{thm:gt} yields a Goldblatt-Thomason theorem
  with respect to this different notion of prime filter extension.

\begin{theorem}
  Let $\ms{K}$ be a class of IM-frames closed
  under $\sigma$-prime filter extensions.
  Then $\ms{K}$ is axiomatic iff it reflects $\sigma$-prime filter
  extensions and is closed under disjoint unions, generated subframes and
  p-morphic images.
\end{theorem}

%--------------------------------------------------------------------------------
\subsection{Non-Normal Intuitionistic Modal Logic}

  Neighbourhood semantics
  is used to accommodate for non-normal modal
  operators~\cite{Sco70,Mon70,Che80,Pac17}.
  Dalmonte, Grellois and Olivett recently put forward an intuitionistic
  analogue~\cite{DalGreOli20}
  to interpret the extension of
  intuitionistic logic with unary modalities $\Box$ and $\Diamond$ which
  a priori do not satisfy any interaction axioms.

  The ordered sets underlying the neighbourhood semantics from~\cite{DalGreOli20}
  are allowed to be preorders. Conforming to our general framework,
  we shall assume them to be posets.
  However, as mentioned in the introduction, we can obtain exactly the same
  (dialgebraic) results when replacing posets with preorders.
  
  We use $\wp$ to denote the (covariant) powerset
  functors on $\cat{Set}$.

\begin{definition}
  A \emph{coupled intuitionistic neighbourhood frame} or \emph{CIN-frame} is a tuple
  $(X, \leq, N_{\Box}, N_{\Diamond})$ such that $(X, \leq)$ is an intuitionistic
  Kripke frame and $N_{\Box}, N_{\Diamond}$ are functions
  $X \to \wp\wp X$ such that for all $x, y \in X$:
  $$
    x \leq y
    \quad\text{implies}\quad N_{\Box}(x) \subseteq N_{\Box}(y)
    \quad\text{and}\quad N_{\Diamond}(x) \supseteq N_{\Diamond}(y).
  $$
  A CIN-morphism
  $f : (X, \leq, N_{\Box}, N_{\Diamond}) \to (X', \leq', N'_{\Box}, N'_{\Diamond})$
  is a p-morphism $f : (X, \leq) \to (X', \leq')$ where
  for all $N \in \{ N_{\Box}, N_{\Diamond} \}$, $x \in X$, $a' \in \power X'$,
  $f^{-1}(a') \in N(x)$ iff $a' \in N'(f(x))$.
  $\cat{CIN}$ denotes the category of CIN-frames and -morphisms.
\end{definition}

  The language $\lan{L}_{\Box\!\Diamond}$ extending the intuitionistic
  language with unary modalities $\Box$ and $\Diamond$ can be
  interpreted in models based on CIN-frames, where 
  \begin{align*}
    x \Vdash \Box\phi &\iff \llb \phi \rrb \in N_{\Box}(x), &
    x \Vdash \Diamond\phi &\iff X \setminus \llb \phi \rrb \notin N_{\Diamond}(x).
  \end{align*}
  We now view this dialgebraically:

\begin{definition}
  Define $\fun{N} : \cat{Krip} \to \cat{Pos}$
  on objects $(X, \leq)$ by 
  $\fun{N}(X, \leq) = (\power\power X, \subseteq) \times (\power\power X, \supseteq)$,
  and on morphisms $f : (X, \leq) \to (X', \leq')$
  by
  $$
    \fun{N}f(W_1, W_2) = \big( \{ a_1' \in \wp X' \mid f^{-1}(a_1') \in W_1 \}, 
                       \{ a_2' \in \wp X' \mid f^{-1}(a_2') \in W_2 \} \big).
  $$
\end{definition}

\begin{theorem}\label{thm:CIN-dialg}
  We have $\cat{CIN} \cong \cat{Dialg}(\fun{i}, \fun{N})$.
\end{theorem}
\begin{proof}
  The isomorphism on objects is obvious. The isomorphism on
  morphisms follows from a computation similar to that in
  the proof of Thm.~\ref{thm:Mon-dialg}.
\end{proof}

  The modal operators $\Box, \Diamond$ are induced by
  $\lambda^{\Box}, \lambda^{\Diamond} : \fun{Up} \circ \fun{i} \to \fun{Up} \circ \fun{N}$,
  where
  \begin{align*}
    \lambda_{(X, \leq)}^{\Box}(a)
      &= \{ (W_1, W_2) \in \fun{N}(X, \leq) \mid a \in W_1 \} \\
    \lambda_{(X, \leq)}^{\Diamond}(a)
      &= \{ (W_1, W_2) \in \fun{N}(X, \leq) \mid X \setminus a \notin W_2 \}
  \end{align*}
  
  Unravelling the definition of a disjoint union of (the dialgebraic
  renderings of) CIN-frames shows that it is computed similar to
  Def.~\ref{def:mif-coprod}.
  Generated subframes and p-morphic images are defined
  by means of CIN-morphisms.

  Since $\Box$ and $\Diamond$ only satisfy the congruence rule,
  the algebraic semantics is given by dialgebras for the functor
  $\fun{L}^{\Box\!\Diamond} : \cat{HA} \to \cat{DL}$ that sends
  $A$ to the free distributive lattice generated by
  $\{ \dbox a, \ddiamond a \mid a \in A \}$.
  The induced natural transformation
  $\rho^{\Box\!\Diamond} : \fun{L}^{\Box\!\Diamond} \circ \upp \to \up \circ \fun{N}$
  is defined on components via
  $\rho^{\Box\!\Diamond}_{(X, \leq)}(\dbox a) = \lambda_{(X, \leq)}^{\Box}(a)$ and
  $\rho^{\Box\!\Diamond}_{(X, \leq)}(\ddiamond a) = \lambda_{(X, \leq)}^{\Diamond}(a)$.
  Akin to Sec.~\ref{subsec:mon-I} we find
  $\dbox a \in (\rho_A^{\Box\!\Diamond})^{\flat}(W_1, W_2)$ iff
  $\theta'_A(a) \in W_1$ and
  $\ddiamond a \in (\rho_A^{\Box\!\Diamond})^{\flat}(W_1, W_2)$
  iff $\fun{pf'}A \setminus \theta'_A(a) \notin W_1$
  for all $A \in \cat{HA}$, $(W_1, W_2) \in \fun{N}(\fun{pf'}A)$ and $a \in A$.

\begin{definition}
  For a Heyting algebra $A$, define
  $$
    \tau_A
      : \fun{pf}(\fun{L}^{\Box\!\Diamond}A) \to \fun{N}(\fun{pf'}A)
      : Q \mapsto \big( \{ \theta'_A(a) \mid \dbox a \in Q \}, \{ \fun{pf'}A \setminus \theta'_A(a) \mid \ddiamond a \notin Q \} \big).
  $$
\end{definition}

  Then $\tau = (\tau_A)_{A \in \cat{HA}}$ defines a natural transformation
  $\fun{pf} \circ \fun{L}^{\Box\!\Diamond} \to \fun{N} \circ \fun{pf'}$.
  It follows from the definitions that
  $(\rho^{\Box\!\Diamond})^{\flat} \circ \tau
   = \fun{id}_{\fun{pf} \circ \fun{L}^{\Box\!\Diamond}}$.
  We get the following definition of $\tau$-prime filter extensions
  and Goldblatt-Thomason theorem.

\begin{definition}
  The $\tau$-prime filter extension of a CIN-frame
  $\mo{X} = (X, \leq, N_{\Box}, N_{\Diamond})$ is given by
  $\pe\mo{X} = (X^{pe}, \subseteq, N_{\Box}^{pe}, N_{\Diamond}^{pe})$,
  where for $\ff{q} \in X^{pe}$ we have
  \begin{align*}
    N_{\Box}^{pe}(\ff{q})
      &= \{ \theta'_{\upp(X, \leq)}(a) \in \wp X^{pe}
             \mid a \in \up(X, \leq) \text{ and } \Box_N(a) \in \ff{q} \} \\
    N_{\Diamond}^{pe}(\ff{q})
      &= \{ X^{pe} \setminus \theta'_{\upp(X, \leq)}(a) \in \wp X^{pe}
             \mid a \in \up(X, \leq) \text{ and } \Diamond_N(a) \in \ff{q} \}
  \end{align*}
  Here $\Box_N(a) = \{ x \in X \mid a \in N_{\Box}(x) \}$
  and $\Diamond_N(a) = \{ x \in X \mid X \setminus a \notin N_{\Diamond}(x) \}$.
\end{definition}

\begin{theorem}
  Let $\ms{K}$ be a class of CIN-frames closed
  under $\tau$-prime filter extensions.
  Then $\ms{K}$ is axiomatic iff it reflects $\tau$-prime filter
  extensions and is closed under disjoint unions, generated subframes and
  p-morphic images.
\end{theorem}

%--------------------------------------------------------------------------------
\subsection{Heyting-Lewis Logic}

  Finally we discuss Heyting-Lewis logic, the extension of intuitionistic
  logic with a binary strict implication operator
  $\sto$~\cite{LitVis18,LitVis19,GroLitPat21}.

\begin{definition}
  A \emph{strict implication frame} is a tuple $(X, \leq, R_s)$,
  where $(X, \leq)$ is an intuitionistic
  Kripke frame and $R_s$ is a relation on $X$ such that $x \leq y R_s z$ implies
  $x R_s z$.
  Morphisms between them are functions that are p-morphisms with respect to 
  both orders.
  Models are defined as expected, and $\sto$ is interpreted via
  $$
    x \Vdash \phi \sto \psi \iff \text{for all } y \in X,
      \text{ if }xR_sy \text{ and } y \Vdash \phi \text{ then } y \Vdash \psi.
  $$
\end{definition}

  Strict implication frames can be modelled as
  $(\fun{i}, \fun{P}_s)$-dialgebras, where
  $\fun{P}_s : \cat{Krip} \to \cat{Pos}$ is the functor that sends
  $(X, \leq)$ to $(\power X, \subseteq)$ ($\power$ denotes the covariant
  powerset functor)
  and a p-morphism $f$ to $\power f$.
  The modality $\sto$ can then be defined via the binary predicate lifting
  $\lambda^{\sto}$, given on components by
  $$
    \lambda^{\sto}_{(X, \leq)}(a, b)
      = \{ c \in \fun{P_s}(X, \leq) \mid c \cap a \subseteq b \}.
  $$
  Disjoint unions, generated subframes and p-morphic images are
  defined as for $\Box$-frames.
  
  The algebraic semantics for this logic given in~\cite[Def.~III.1]{GroLitPat21}
  can be modelled dialgebraically in a similar way as we have seen above.
  Computation of the natural transformation $\rho^{\sto}$ is, by now, routine.
  Examining the proof of the duality for Heyting-Lewis logic sketched
  in~\cite[Section III-D]{GroLitPat21}, we can compute
  a one-sided inverse $\tau$ to $(\rho^{\sto})^{\flat}$.
  We suppress the details, but do give the resulting notion of prime filter
  extension:

\begin{definition}
  The \emph{prime filter extension} of a strict implication frame $(X, \leq, R_s)$
  is given by the frame $(X^{pe}, \subseteq, R_s^{pe})$, with $R_s^{pe}$ defined by
  $$
    \fil{p} R_s^{pe} \fil{q} \iff \forall a, b \in \up(X, \leq),
      \text{if } a \sto_R b \in \fil{p} \text{ and } a \in \fil{q}
      \text{ then } b \in \fil{q}
  $$
  where $a \sto_R b = \{ x \in X \mid R[x] \cap a \subseteq b \}$.
\end{definition}

  With this notion of prime filter extension, Thm.~\ref{thm:gt} instantiates to:
  
\begin{theorem}
  A class $\ms{K}$ of strict implication frames that is closed under prime
  filter extensions is axiomatic iff it reflects
  prime filter extensions and is closed under disjoint unions, generated
  subframes and p-morphic images.
\end{theorem}

%================================================================================
\section{Conclusions}

  We have given a general way to obtain Goldblatt-Thomason theorems for
  modal intuitionistic logics, using the framework of dialgebraic logic.
  Subsequently, we applied the general result to several concrete modal
  intuitionistic logics.
  The results in this paper can be generalised in several directions.

\begin{description}
  \item[More applications.]
        The general Goldblatt-Thomason theorem can also be instantiated to
        $\Diamond$-frames and $\Box\Diamond$-frames~\cite{WolZak99}.
        Using preorders instead of posets, we can obtain Goldblatt-Thomason
        theorems for ((strictly) condensed) $H\Box$-frames and
        $H\Box\Diamond$ frames used by Bo\v{z}i\'{c} and
        Do\v{s}en~\cite{BozDos84}.
  \item[More base logics.]
        The framework of dialgebraic logic is not restricted to an intuitionistic
        base. Generalising the results from this paper, we can obtain a general
        Goldblatt-Thomason theorem that also covers modal bi-intuitionistic
        logics~\cite{GroPat19} and modal lattice logics~\cite{BezEA22}.
        Moreover, this would also cover coalgebraic logics over a classical
        and a positive propositional base.
        The results in this paper can be generalised to dialgebraic logics
        for different base logics. This would give rise to Goldblatt-Thomason
  \item[Other modal intuitionistic logics]
        The results in the paper do not apply to the modal
        intuitionistic logics investigated by~Fischer Servi~\cite{Fis81},
        Plotkin and Sterling~\cite{PloSti86}, and Simpson~\cite{Sim94},
        because these formalisms are not covered by the dialgebraic approach.
        It would be interesting to see if similar techniques can be applied
        to these logics to still prove Goldblatt-Thomason theorems.
\end{description}
  
\bigskip\noindent
\textit{Acknowledgements.} \,
  I am grateful to the anonymous reviewers for many constructive and helpful comments.

%% Bibliography
{\footnotesize
\bibliographystyle{plainnat}
\bibliography{../../../../biblio.bib}
}

\clearpage
\appendix
\section{Omitted proofs}

  We use the following lemma in the proof of Prop.~\ref{prop:pf-key}.

\begin{lemma}\label{lem:tau-nice}
  Let $\tau$ be a natural transformation such that
  $\rho^{\flat} \circ \tau = \fun{id}_{\fun{pf} \circ \fun{L}}$,
  and $\amo{A} = (A, \alpha) \in \cat{Dialg}(\fun{L}, \fun{j})$.
  Then $\theta'_A : A \to \upp(\fun{pf'}H)$ defines a
  $(\fun{L}, \fun{j})$-dialgebra morphism from $\amo{A}$ to
  $(\amo{A}_{\tau})^+$.
\end{lemma}
\begin{proof}
  This is similar to \cite[Theorem~6.4(1)]{KurRos12}.
  We repeat the argument here.

  Let $\amo{A} = (A, \alpha)$ be a $(\fun{L}, \fun{j})$-dialgebra.
  Then $(\amo{A}_{\tau})^+$ is given by the composition
  $$
  \footnotesize{
    \begin{tikzcd}
      \fun{L}(\upp(\fun{pf'}A))
            \arrow[r, "\rho_{\fun{pf'}A}"]
        & \up(\fun{T}(\fun{pf'}A))
            \arrow[r, "\up\tau_{A}"]
        & \up(\fun{pf}(\fun{L}A))
            \arrow[r, "\up \circ \fun{pf} \alpha"]
        & \up(\fun{pf}(\fun{j}A))
            \arrow[r, equal]
        & [-1em]\fun{j}(\upp(\fun{pf'}A))
    \end{tikzcd}
  }
  $$
  In order to show that $\theta'_A$ is a morphism from
  $\amo{A}$ to $(\amo{A}_{\tau})^+$ we need to show that the
  outer shell of the following diagram commutes:
  $$
  \footnotesize{
  \begin{tikzcd}[column sep=2em]
    \fun{L}A        \arrow[d, "\fun{L}\theta'_A" left]
                    \arrow[rrrr, "\alpha"]
                    \arrow[drr, "\theta_{\fun{L}A}"]
      &&&& [-.5em] \fun{j}A \arrow[d, "\fun{j}\theta'_A"]
                    \arrow[dl, "\theta_{\fun{j}A}" swap] \\
    \fun{L}(\upp(\fun{pf'}A))
          \arrow[r, "\rho_{\fun{pf'}A}" below]
      & \up(\fun{T}(\fun{pf'}A))
          \arrow[r, "\up\tau_A" below]
      & \up(\fun{pf}(\fun{L}A))
          \arrow[r, "\up(\fun{pf}\alpha)" below]
      & \up(\fun{pf}(\fun{j}A))  \arrow[r, equal]
      & \fun{j}(\upp(\fun{pf'}A))
  \end{tikzcd}
  }
  $$
  The right triangle commutes by definition.
  The middle square commutes by naturality of $\theta$.
  So we are left to prove that
  $\theta_{\fun{L}A} = \up\tau_A \circ \rho_{\fun{pf'}A} \circ \fun{L}\theta'_A$.
  
  Since $\rho^{\flat} \circ \tau = \id$,
  hence $\up\tau \circ \up\rho^{\flat} = \id_{\up}$,
  it suffices to prove that
  $\up\rho^{\flat} \circ \theta_{\fun{L}A} = \rho_{\fun{pf'}A} \circ \fun{L}\theta'_A$.
  (The result then follows from composing both sides with $\up\tau_A$ on the left.)
  This is precisely the outer shell of the diagram
  $$
  \footnotesize{
  \begin{tikzcd}
    \fun{L}A
          \arrow[r, "\fun{L}\theta_A'"] \arrow[d, "\theta_{\fun{L}A}" left]
      & \fun{L}(\upp(\fun{pf'}A))
          \arrow[r, "\rho_{\fun{pf'}A}"]
          \arrow[d, "\theta_{\fun{L}(\upp(\fun{pf'}A))}" left]
      & \up(\fun{T}(\fun{pf'}A))
          \arrow[dr, bend left=20, "\id"]
          \arrow[d, "\theta_{\up(\fun{T}(\fun{pf'}A))}" left]
      & \\
    \up(\fun{pf}(\fun{L}H))
          \arrow[r, "\up(\fun{pf}(\fun{L}\theta_A'))"]
          \arrow[rrr, bend right=12, "\up\rho_A^{\flat}" above]
      & \up(\fun{pf}(\fun{L}(\upp(\fun{pf'}A))))
          \arrow[r, "\up(\fun{pf}\rho_{\fun{pf'}A})"]
      & \up(\fun{pf}(\up(\fun{T}(\fun{pf'}A))))
          \arrow[r, "\up\eta_{\fun{T}(\fun{pf'}A)}"]
      & \up(\fun{T}(\fun{pf'}A))
  \end{tikzcd}
  }
  $$
  Here the bottom square commutes by definition of $\rho^{\flat}$.
  The other two squares commute by naturality of $\theta$ and the triangle
  on the right commutes because $\theta$ and $\eta$ are the units of a
  dual adjunction.
\end{proof}

\begin{proof}[Proof of Proposition~\ref{prop:pf-key}.]
  Recall that $\theta'_{\upp(X, \leq)}(\llb \phi \rrb^{\mo{M}}) =
  \{ \ff{p} \in \fun{pf'}(\upp(X, \leq)) \mid \llb \phi \rrb^{\mo{M}} \in \ff{p} \}$.
  So the first item is equivalent to
  $$
    \llb \phi \rrb^{\pe\mo{M}} = \theta'_{\upp(X, \leq)}(\llb \phi \rrb^{\mo{M}}),
  $$
  where we view truth sets of formulae as elements in the relevant complex algebras
  (cf.~Prop.~\ref{prop:complex-alg}).
  The proof proceeds by induction on the structure of $\phi$.
  If $\phi = q \in \Prop$ then the statement holds by definition of
  $V^{pe}$. The cases $\phi = \top$ and $\phi = \bot$ hold by definition
  of a prime filter.
  
  If $\phi$ is of the form $\phi_1 \star \phi_2$, where
  $\star \in \{ \wedge, \vee, \to \}$ then we use Lem.~A to find
  \begin{align*}
    \llb \phi_1 \star \phi_2 \rrb^{\pe\mo{M}}
      &= \llb \phi_1 \rrb^{\pe\mo{M}} \star \llb \phi_1 \rrb^{\pe\mo{M}}
          \\
      &= \theta'_{\upp(X, \leq)}(\llb \phi_1 \rrb^{\mo{M}}) \star
         \theta'_{\upp(X, \leq)}(\llb \phi_2 \rrb^{\mo{M}})
         &\text{(IH)} \\
      &= \theta'_{\upp(X, \leq)}(\llb \phi_1 \rrb^{\mo{M}} \star
         \llb \phi_2 \rrb^{\mo{M}}) \\
      &= \theta'_{\upp(X, \leq)}(\llb \phi_1 \star \phi_2 \rrb^{\mo{M}})
  \end{align*}
  The case where $\phi = \heartsuit^{\lambda}(\phi_1, \ldots, \phi_n)$ follows
  from a similar computation, using the fact that $\theta'_{\up(X, \leq)}$
  preserves operators of the form $\dheartsuit^{\lambda}$.
  
  Item (ii) follows from Item (i) and the definition of
  $\eta_{(X, \leq)}(x)$ via
  $$
    \mo{M}, x \Vdash \phi
      \iff x \in \llb \phi \rrb^{\mo{M}}
      \iff \llb \phi \rrb^{\mo{M}} \in \eta_{(X, \leq)}(x)
      \iff \pe\mo{M}, \eta_{(X, \leq)}(x) \Vdash \phi.
  $$
  For Item (iii),
  let $V$ be any valuation for $\mo{X}$ and $x \in X$.
  By assumption $(\pe\mo{X}, V^{pe}), \theta'_{(X, \leq)}(x) \Vdash \phi$,
  so by Item (ii)
  $(\mo{X}, V), x \Vdash \phi$
  and hence $\mo{X} \Vdash \phi$.
\end{proof}

\begin{proof}[Proof of Lemma~\ref{lem:Box-canonical}.]
  Let $A$ be a Heyting algebra.
  Recall that
  $\theta'_A(a) = \{ \ff{q} \in \fun{pf'}A \mid a \in \ff{q} \}$.
  Using this we can rewrite
  $\tau^{\Box}_A : \fun{pf}(\fun{L}^{\Box}A) \to \fun{P_{up}}(\fun{pf'}A)$
  as
  \begin{equation}\label{eq:box-canonical}\tag{$\star$}
    \tau^{\Box}_A(Q) = \bigcap \{ \theta'_A(a) \mid a \in A, \dbox a \in Q \}.
  \end{equation}
  Since $\theta'_A(a)$ is an upset of $\fun{pf'}A$,
  $\tau^{\Box}_A(Q)$ is also an upset of $\fun{pf'}A$, hence in $\fun{P_{up}}(\fun{pf'}A)$.
  The elements of $\fun{pf}(\fun{L}^{\Box}A)$ are ordered by inclusion.
  If $Q, Q' \in \fun{pf}(\fun{L}^{\Box}A)$ and $Q \subseteq Q'$
  then it follows immediately that $\tau^{\Box}_A(Q) \supseteq \tau^{\Box}_A(Q')$.
  Since $\fun{P_{up}}(\fun{pf'}A)$ is ordered by reverse inclusion,
  so $\tau^{\Box}_A$ is a morphism of $\cat{Pos}$.
  
  For naturality,
  let $h : A \to B$ be a Heyting homomorphism.
  We need that
  $$
    \begin{tikzcd}
      \fun{pf}(\fun{L}^{\Box}A)
            \arrow[r, "\tau_A^{\Box}"]
        & \fun{P_{up}}(\fun{pf'}A) \\ [-.25em]
      \fun{pf}(\fun{L}^{\Box}B)
            \arrow[r, "\tau_B^{\Box}"]
            \arrow[u, "(\fun{L}^{\Box}h)^{-1}" left]
        & \fun{P_{up}}(\fun{pf'}B)
            \arrow[u, "\fun{P_{up}}(h^{-1})" right]
    \end{tikzcd}
  $$
  commutes.
  Let $Q \in \fun{pf}(\fun{L}^{\Box}B)$, $\ff{q} \in \fun{pf'}A$,
  and suppose $\ff{q} \in \tau_A^{\Box}(\fun{L}^{\Box}h)^{-1}(Q)$.
  To show $\ff{q} \in \fun{P_{up}}(h^{-1})(\tau_B^{\Box}(Q))$
  it suffices to find a prime filter $\ff{p} \in \tau_B^{\Box}(Q)$
  such that $h^{-1}(\ff{p}) \subseteq \ff{q}$,
  because $\tau_A^{\Box}(\fun{L}^{\Box}h)^{-1}(Q)$ is 
  an upset of $\fun{pf'}A$.
  Define $F = \{ b \in B \mid \dbox b \in Q \}$
  and $I = \{ b \in B \mid \exists a \in A \setminus \ff{q} \text{ s.t. } b \leq c \}$.
  If $I \cap F \neq \emptyset$ then there exists $b \in B$ and $a \in A \setminus \ff{q}$
  such that $b \leq h(a)$. Since $\dbox b \in Q$ this implies
  $\dbox h(a) \in Q$ and hence $\dbox a \in (\fun{L}^{\Box}h)^{-1}(Q)$.
  But then $a \in \ff{q}$ because $\ff{q} \in \tau^{\Box}_A((\fun{L}^{\Box}h)^{-1}(Q))$,
  a contradiction.
  So $I \cap F = \emptyset$.
  The prime filter lemma then gives a prime filter $\ff{p}$
  containing $F$ and disjoint from $I$.
  This satisfies $\ff{p} \in \tau_B^{\Box}(Q)$
  and $h^{-1}(\ff{p}) \subseteq \ff{q}$ by design.
  
  Conversely, suppose $\ff{q} \in \fun{P_{up}}(h^{-1})(\tau_B^{\Box}(Q))$.
  Then there exists a $\ff{p} \in \tau_B^{\Box}(Q)$ such that
  $\ff{q} = h^{-1}(\ff{p})$.
  We show that $\ff{q} \in \tau_A^{\Box}(\fun{L}^{\Box}h)^{-1}(Q)$.
  Let $a \in A$ and suppose $\dbox a \in (\fun{L}^{\Box}h)^{-1}(Q)$.
  Then $\dbox h(a) = \fun{L}^{\Box}h(\dbox a) \in Q$
  so $h(a) \in \ff{p}$. But this implies $a \in \ff{q} = h^{-1}(\ff{p})$.
  So by definition $\ff{q} \in \tau_A^{\Box}(\fun{L}^{\Box}h)^{-1}(Q)$.
  
  Finally, we show that
  $(\rho^{\Box}_A)^{\flat} \circ \tau^{\Box}_A
    = \fun{id}_{\fun{pf} \circ \fun{L}^{\Box}A}$
  for $A \in \cat{HA}$. Let $Q \in \fun{pf}(\fun{L}^{\Box}A)$.
  Since elements of $\fun{pf}(\fun{L}^{\Box}A)$ are determined
  uniquely by the generators of the form $\dbox a$ they contain,
  it suffices to show that
  $\dbox a \in Q$ iff $\dbox a \in (\rho^{\Box}_A)^{\flat}(\tau^{\Box}_A(Q))$.
  Because of the computation in Exm.~\ref{exm:box-frame-tau} this
  is equivalent to showing
  $\dbox a \in Q$ iff $\tau^{\Box}_A(Q) \subseteq \theta'_A(a)$.
  The direction from left to right follows from \eqref{eq:box-canonical}.
  For the converse, suppose $\dbox a \notin Q$.
  Let $F = \{ b \in A \mid \Box b \in Q \}$ and $I = \{ c \in A \mid c \leq a \}$.
  Then $F$ is a filter and $I$ is an ideal of $A$, and $F \cap I = \emptyset$.
  By the prime filter lemma we obtain some $\ff{q} \in \fun{pf'}A$ extending
  $F$ and disjoint from $I$. This implies that $\ff{q} \in \tau_A^{\Box}(Q)$ while
  $\ff{q} \notin \theta'_A(a)$,
  so that $\tau^{\Box}_A(Q) \not\subseteq \theta'_A(a)$.
\end{proof}

\begin{proof}[Proof of Lemma~\ref{lem:mon-tau-natural}.]
  Throughout this proof use the fact that $\fun{pf'}A$ forms
  an Esakia space (which in particular is a Stone space), with a topology
  generated by sets of the form $\theta'_A(a)$ and their complements
  \cite[Sec.~2.3.3]{Bez06}.
  Furthermore, we note that for any Heyting homomorphism $h : A \to B$
  we have
  \begin{equation}\label{eq:h-theta}\tag{$\dagger$}
    \theta'_B(h(a)) = (h^{-1})^{-1}(\theta_A'(a))
  \end{equation}
  
  We first prove naturality of $\tau^{\vartri}$. Let $h : A \to B$
  be a Heyting homomorphism. We need to show that
  the following diagram commutes:
  $$
    \begin{tikzcd}
      \fun{pf}(\fun{L}^{\vartri}A)
            \arrow[r, "\tau_A^{\vartri}"]
        & \fun{P_{up}}(\fun{pf'}A) \\ [-.25em]
      \fun{pf}(\fun{L}^{\vartri}B)
            \arrow[r, "\tau_B^{\vartri}"]
            \arrow[u, "(\fun{L}^{\vartri}h)^{-1}" left]
        & \fun{P_{up}}(\fun{pf'}B)
            \arrow[u, "\fun{P_{up}}(h^{-1})" right]
    \end{tikzcd}
  $$
  Let $Q \in \fun{pf}(\fun{L}^{\vartri}B)$ and
  $D \in \fun{Up}(\fun{pf'}A)$. We go by the items of Def.~\ref{def:mon-tau}.
  \begin{itemize}
    \item If $D = \theta'_A(a)$ for some $a \in A$ then
          \begin{align*}
            \theta'_A(a) \in \tau^{\vartri}_A(\fun{L}^{\vartri}h)^{-1}(Q)
              &\iff \dtriangle a \in (\fun{L}^{\vartri}h)^{-1}(Q)
                    &\text{(Def.~$\tau^{\vartri}$)} \\
              &\iff (\fun{L}^{\vartri}h)(\dtriangle a) \in Q \\
              &\iff \dtriangle h(a) \in Q
                    &\text{(Def.~of $\fun{L}^{\vartri}$)} \\
              &\iff \theta_B'(h(a)) \in \tau_B^{\vartri}(Q)
                    &\text{(Def.~of $\tau^{\vartri}$)} \\
              &\iff (h^{-1})^{-1}(\theta'_A(a)) \in \tau_B^{\vartri}(Q)
                    &\text{(By~\eqref{eq:h-theta})} \\
              &\iff \theta'_A(a) \in \fun{M}(h^{-1})(\tau_B^{\vartri}(Q))
                    &\text{(Def.~of $\fun{M}$)}
          \end{align*}
    \item Suppose $D$ is closed in $\fun{pf'}A$.
          If $D \in \tau_A^{\vartri}(\fun{L}^{\vartri}h)^{-1}(Q)$,
          then for all $a \in A$, $D \subseteq \theta'_A(a)$ implies
          $\dtriangle a \in (\fun{L}^{\vartri}h)^{-1}(Q)$, i.e.~$\dtriangle h(a) \in Q$.
          In order to prove that $D \in \fun{M}(h^{-1})(\tau_B^{\vartri}(Q))$,
          we need to show that
          $(h^{-1})^{-1}(D) \in \tau_B^{\vartri}$.
          Since $h^{-1}$ is an Esakia morphism (hence continuous),
          $(h^{-1})^{-1}(D)$ is closed in $\fun{pf'}B$, so it suffices to show
          that for all $b \in B$, $(h^{-1})^{-1}(D) \subseteq \theta'_B(b)$ implies
          $\dtriangle b \in Q$.
          Let $b \in B$ be such that $(h^{-1})^{-1}(D) \subseteq \theta'_B(b)$.
          Then since $D$ is closed we have
          $$
            \bigcap \big\{ (h^{-1})^{-1}(\theta'_A(a))
              \mid a \in A, D \subseteq \theta'_A(a) \}
              \subseteq \theta'_B(b) \big\}.
          $$
          Using~\eqref{eq:h-theta} and compactness of $\fun{pf'}B$ 
          we can find 
          $a_1, \ldots, a_n \in A$ such that
          $$
            \theta'_B(h(a_1 \wedge \cdots \wedge a_n))
              = \theta'_B(h(a_1)) \cap \cdots \cap \theta'_B(h(a_n))
              \subseteq \theta'_B(b).
          $$
          As a consequence of Esakia duality it follows that
          $h(a_1 \wedge \cdots \wedge a_n) \leq b$.
          Since $D \subseteq \theta'_A(a_1 \wedge \cdots \wedge a_n)$,
          we have
          $\dtriangle (a_1 \wedge \cdots \wedge a_n) \in (\fun{L}^{\vartri}h)^{-1}(Q)$,
          so $\dtriangle(h(a_1 \wedge \cdots \wedge a_n)) \in Q$.
          Monotonicity of $\dtriangle$ now implies $\dtriangle b \in Q$.
          
          Conversely, if $D \in \fun{M}(h^{-1})(\tau_B^{\vartri}(Q))$ then
          a similar but easier argument shows that 
          $D \in \tau_A^{\vartri}(\fun{L}^{\vartri}h)^{-1}(Q)$.
    \item Finally, suppose $D$ is any upset.
          If $D \in \tau_A^{\vartri}(\fun{L}^{\vartri}h)^{-1}(Q)$
          then there exists a closed upset $C$ such that
          $C \subseteq D$ and $C \in \tau_A^{\vartri}(\fun{L}^{\vartri}h)^{-1}(Q)$.
          This implies $C \in \fun{M}(h^{-1})(\tau_B^{\vartri}(Q))$,
          so that $(h^{-1})^{-1}(C) \in \tau_B^{\vartri}(Q)$.
          Since $(h^{-1})^{-1}(C)$ is closed again and
          $(h^{-1})^{-1}(C) \subseteq (h^{-1})^{-1}(D)$ we have
          $(h^{-1})^{-1}(D) \in \tau_B^{\vartri}(Q)$,
          and therefore $D \in \fun{M}(h^{-1})(\tau_B^{\vartri}(Q))$.
          
          Conversely, suppose
          $D \in \fun{M}(h^{-1})(\tau_B^{\vartri}(Q))$.
          Then there exists a closed upset $C \in \tau_B^{\vartri}(Q)$
          such that $C \subseteq (h^{-1})^{-1}(D)$.
          Define $C' = h^{-1}[C]$ to be direct image of $C$ under $h^{-1}$.
          Since $h^{-1}$ is an Esakia morphism it sends closed upsets
          to closed upsets. Furthermore $C \subseteq (h^{-1})^{-1}(C')$ so
          $C' \in \fun{M}(h^{-1})(\tau_B^{\vartri}(Q))$.
          This implies $C' \in \tau_A^{\vartri}(\fun{L}^{\vartri}h)^{-1}(Q)$.
          By design $C' \subseteq D$, hence
          $D \in \tau_A^{\vartri}(\fun{L}^{\vartri}h)^{-1}(Q)$.
  \end{itemize}

  Next we prove that $(\rho^{\vartri})^{\flat}_A \circ \tau_A = \fun{id}_{\fun{pf}(\fun{L}^{\vartri}A)}$ for $A \in \cat{HA}$.
  It follows from the definitions of $\rho^{\flat}$ and $\tau$ that for any
  Heyting algebra $A$, $a \in A$ and prime filter $Q \in \fun{pf}(\fun{L}^{\vartri}A)$
  we have
  $\theta_A'(a) \in \rho^{\flat}_A(\tau_A(Q))$
  iff $\dtriangle a \in \tau_A(Q)$ iff $\theta_A'(a) \in Q$.
  Since elements of $\fun{pf}(\fun{L}^{\vartri}A)$ are determined
  uniquely by the elements of the form $\dtriangle a$ they contain,
  this proves the lemma.
\end{proof}

\end{document}